\documentclass{article}
\usepackage{color}
\usepackage{amsmath}
\usepackage{amssymb}
\usepackage{bm}
\usepackage{upgreek}
\usepackage{multirow}
\usepackage{stackengine}
\usepackage{comment}
\usepackage{multicol,caption}
\usepackage[a4paper, total={170mm,257mm}, left=20mm, top=20mm,]{geometry} 
\usepackage{graphicx}
\newenvironment{Figure}
  {\par\medskip\noindent\minipage{\linewidth}}
  {\endminipage\par\medskip}

\newcommand{\numcircle}[1]{\mathbin{\stackMath{\stackinset{c}{0pt}{c}{0pt}{#1}{\bigcirc}}}}

\def\ar{\begin{array}}
	\def\arr{\end{array}}
\def\be{\begin{eqnarray}}
	\def\en{\end{eqnarray}}
\def\bee{\begin{equation}}
	\def\ee{\end{equation}}
\graphicspath{{./Figs/}}

\title{FFT-acceleration and stabilization of the 3D Marching-on-in-Time Contrast Current Density Volume Integral Equation for scattering from high contrast dielectrics}

\author{Petrus~W.N.~van~Diepen \thanks{Corresponding author:~Petrus~Wilhelmus~Nicolaas~van~Diepen (p.w.n.v.diepen@tue.nl).} \thanks{Department of Electrical Engineering, Eindhoven University of Technology, The Netherlands.}, Martijn~C.~van~Beurden\footnotemark[2] and Roeland~J.~Dilz\footnotemark[2]}

\begin{document}
\maketitle

\begin{abstract}
An implicit causal space-time Galerkin scheme applied to the contrast current density volume integral equation gives rise to a marching-on-in-time scheme known as the MOT-JVIE, which is accelerated and stabilized via a fully embedded FIR filter to compute the electromagnetic scattering from high permittivity dielectric objects discretized with over a million voxels. A review of two different acceleration approaches previously developed for two-dimensional time-domain surface integral equations based on fast Fourier transforms (FFTs), leads to an understanding why these schemes obtain the same order of acceleration and the extension of this FFT-acceleration to the three-dimensional MOT-JVIE. The positive definite stability analysis (PDSA) for the MOT-JVIE shows that the number of voxels for a stable MOT-JVIE discretization is restricted by the finite precision of the matrix elements. The application of the PDSA provides the insight that stability can be enforced through regularization, at the cost of accuracy. To minimize the impact in accuracy, FIR-regularization is introduced, which is based on low group-delay linear-phase high-pass FIR-filters. We demonstrate the capabilities of the FFT-accelerated FIR-regularized MOT-JVIE for a number of numerical experiments with high permittivity dielectric scatterers.
\end{abstract}

\begin{multicols}{2}
\section{Introduction}\label{section label}
Time-domain Maxwell solvers are most suitable for full-wave electromagnetic simulations dealing with either short-time transient analysis, ultra-wideband excitations, time-modulated materials, non-linearity in the material properties, multiphysics or a combination of these ~\cite{Ren2022,Miller1999}. The differential-equation-based time domain (DETD) Maxwell solvers are preferred for these types of simulations~\cite{Jin2019,Sankaran2019,Ren2022}, in particular the finite difference time domain (FDTD) method, finite element time domain (FETD) method and discontinuous Galerkin time domain (DGTD) method. Time domain integral equations (TDIEs) are, unlike the differential-equation-based methods, based on the time-domain Green function. Consequently, there are certain advantages of TDIEs over the DETD methods~\cite{Miller1999,Weile2019,Ren2022}: 1) the solution inherently satisfies the radiation condition and no numerical truncation of the boundary is required; 2) the background medium is not included in the computational domain. Together, these advantages result in a reduction of the computational domain as compared to DETD methods, which is especially relevant for simulations with low-frequency content as the included background volume increases in DETD methods~\cite{Berenger2007}.

There are three classes of discretizations for TDIEs~\cite{Ren2022}: convolutional quadrature (CQ)~\cite{Wang2011,Ding2016}, marching-on-in-degree (MOD)~\cite{Chung2004,Shi2011a} and space-time Galerkin~\cite{Rao1991,Gres2001}. The difference between these classes is how they handle the time aspect of the TDIEs to find the solution at any time instance. The CQ schemes are obtained by mapping the Laplace-domain equivalent of the TDIEs to the $\mathcal{Z}$-domain using an implicit Runge-Kutta method and  numerically computing the inverse $\mathcal{Z}$-transform to find the unknown through a marching-on-in-time (MOT) scheme~\cite{Wang2011,Ding2016}, i.e. finding the unknown at increasing time instances. The MOD schemes employ a set of globally defined, i.e. along the entire time axis, orthonormal basis functions to expand and test time~\cite{Chung2004,Shi2011a}, where the unknown at any time instance is collectively represented by all these basis functions. The space-time Galerkin schemes employ locally defined temporal basis and test functions to expand and test in time to obtain a MOT-scheme similar to CQ schemes, but the unknown at any time instance is represented by a limited number of basis functions~\cite{Rao1991,Gres2001}. The advantage of space-time Galerkin schemes over the other two is that these schemes allow for non-linearity in the simulation, unlike the MOD schemes, and are immune to numerical dispersion, unlike CQ schemes~\cite{Ren2022}. 

The space-time Galerkin schemes can be distinguished by how they expand and test time with locally defined temporal functions and we differentiate  implicit causal, implicit non-causal, and explicit schemes.
The implicit causal schemes expand and test in time such that the unknowns in the matrix equation do not depend on future solutions~\cite{Dodson1998,Gres2001}. The solution at the current time step can be found by evaluating the field contributions of solutions at earlier time steps and by subsequently solving a linear system per time step. The sparsity of the matrix to be inverted is proportional to the length of the time step size. The implicit non-causal schemes expand and test in time such that the solution at each time step does depend on the future solutions~\cite{Weile2004,Sayed2015}. As a result, the implicit non-causal schemes require an additional predictor-corrector step. The traditional explicit schemes~\cite{Rao1991,Al-Jarro2012} try to circumvent the step of solving a linear system at every time step in implicit schemes by reducing the sparsity of the pertaining matrix to a diagonal. This comes at the cost of a Courant-Friedrich-Lewy (CFL) like condition, whereas for implicit schemes the smallest time step is governed by the maximum frequency of the excitation~\cite{Ren2022}. A more recent class of explicit schemes has been developed in~\cite{Sayed2020}, which casts the matrix equation in an ordinary differential equation (ODE) employing Lagrange temporal interpolation and the ODE is solved with a predictor-corrector in combination with solving the Gram matrix equation at each prediction and correction step. This explicit scheme has no CFL condition, unlike the traditional explicit methods. 

The space-time Galerkin schemes for TDIEs can be applied to a surface integral formulation (TDSIE) and to a volume integral formulation (TDVIE)~\cite{Jin2010}. The TDSIEs are used to represent wave propagation in homogeneous media whereas TDVIEs can also handle wave propagation in inhomogeneous media. The latter is important if time-varying or field dependent material properties are to be considered, where inhomogeneity is common. TDVIEs can be applied to polarized, magnetized, and lossy propagation media, however research in space-time Galerkin schemes based on TDVIEs has so far focused on polarizing media. 

A long-standing problem of space-time Galerkin schemes based on TDVIEs has been to maintain stability for an increasing dielectric contrast, i.e. the permittivity of a dielectric object is large compared to that of background medium~\cite{Sayed2015}. Three schemes have shown promising results: an implicit non-causal scheme~\cite{Sayed2015}, an explicit scheme~\cite{Sayed2020} and an implicit causal scheme~\cite{VanDiepen2024}. The advantage of the latter over the other two is that the method does not require a predictor-corrector method to find the solution at each time step. 

The implicit causal space-time Galerkin scheme in~\cite{VanDiepen2024} is referred to as the marching-on-in-time contrast current density volume integral equation (MOT-JVIE). Like most MOT schemes, the computation speed is limited by the linear algebra operations required to compute the field produced by earlier computed solutions. Two popular methods that achieve acceleration of these operations exist in literature, i.e. the (multi-level) plane-wave time domain (PWTD) method~\cite{Ergin1999, Shanker2003, Shanker2004} and the FFT-acceleration~\cite{Hairer1985,Yilmaz2001,Yilmaz2002,Yilmaz2002a,Yilmaz2004}. Although both techniques can be applied to the MOT-JVIE, the FFT-acceleration is the most straightforward to incorporate. The MOT-JVIE applies piece-wise constant basis functions defined on voxels~\cite{VanDiepen2024}, which will restrict the spatial part of the contrast current density to the appropriate solution space, i.e. $L^2(\mathbb{R}^3)$~\cite{Beurden2003,Beurden2007,VanBeurden2008}. The voxels can be defined on a regularized grid, which introduces a discrete spatial shift invariance in the matrices, which is necessary for FFT-acceleration~\cite{Yilmaz2001}.

Two implementations of FFT-acceleration for TDSIEs exist, the one presented in~\cite{Yilmaz2002} and the one presented in~\cite{Yilmaz2002a}. At first glance, these FFT-accelerations seem completely different. However, they achieve the same order of acceleration. To understand their differences and why they achieve the same acceleration, we review these techniques, before we apply them to the MOT-JVIE. Up to the authors knowledge, neither technique has been applied yet to space-time Galerkin schemes based on TDVIEs.

The FFT-acceleration of the MOT-JVIE enables simulations for discretizations with a significantly larger number of voxels. In subsequent numerical experiments, an instability appears when the number of voxels is increased. The topic of instability in TDIEs due to an increase in spatial elements is rare, but in~\cite{Shanker2009} it is mentioned that unpublished studies show a similar phenomenon for the M{\"u}ller formulation of the TDSIEs. In the same work, an improvement in stability is obtained by putting more effort in the accurate evaluation of the integrals, which according to~\cite{Shanker2009} illustrates a link between stability and finite precision in the matrix elements. The link between accuracy and stability is further established in~\cite{Wout2013} for TDSIEs. By means of the positive definite stability analysis (PDSA)~\cite{VanDiepen2024a}, we investigate whether the finite precision in the matrix elements of the MOT-JVIE interaction matrices explains the loss of stability for an increasing number of voxels. Some strategies to improve the accuracy were already included in the MOT-JVIE~\cite{VanDiepen2024}, i.e. by choosing a time step larger than the spatial step, which increases the smoothness of the integrand~\cite{Dodson1998}, and by exact evaluation of the radiated fields~\cite{Shanker2009,Wout2013}. Here, we follow a different strategy based on the insights of the PDSA in the form of a regularization to enforce a stability of the MOT-JVIE. 

The paper is organized as follows. After introducing the MOT-JVIE in Section~\ref{sc:MOT-JVIE}, we focus on two parts. First, we review existing FFT-acceleration of the TDSIEs and extend this to the MOT-JVIE in Sections~\ref{sc:SpatialAcceleration} and~\ref{sc:TemporalAcceleration}. Second, to enable stable computation of long time sequences for large bodies, we study the link between finite precision and the number of voxels in the MOT-JVIE and discuss several regularization options to enforce stability and show their impact on the accuracy of the solution in Section~\ref{sc:Stabilization}. The capabilities of the regularized MOT-JVIE are then demonstrated in Section~\ref{sc:AccRegMOTJVIE}. Finally, we draw conclusions in Section~\ref{sc:Conclusion}.

\section{MOT-JVIE} \label{sc:MOT-JVIE}
\subsection{Formulation} \label{sc:Formulation}
A dielectric object in a homogeneous medium, with permittivity $\varepsilon_0$ and permeability $\mu_0$ and resulting wave speed
$c_0 = 1/\sqrt{\varepsilon_0 
\mu_0}$, occupies a volume  $V_\varepsilon$. The permittivity $\varepsilon(\mathbf{r})$ inside the volume is position dependent and is defined relative to that of the background medium, i.e. $\varepsilon = \varepsilon_r(\mathbf{r})\varepsilon_0$ with relative permittivity $\varepsilon_r(\mathbf{r}) \geq 1$. The contrast current density $\mathbf{J}_\varepsilon(\mathbf{r},t)$ inside this volume is then induced by an incident electric field $\mathbf{E}^i(\mathbf{r},t)$ that arrives at the object after $t = 0$. Subsequently, the contrast current density generates the scattered electromagnetic field in the background medium in accordance to the convolution with the Green function~\cite{Jackson1962} defined as
\begin{equation} \label{eq:GreenFunction}
    G(\mathbf{r}',\mathbf{r},t',t) = \frac{\delta(t - \tau)}{4 \pi R},
\end{equation}
with retarded time $\tau = t-\frac{R}{c_0}$ and distance $R = |\mathbf{r}-\mathbf{r}'|$. The scattered magnetic field strength is then represented by 
\begin{equation} \label{eq:ScatteredMagneticFieldStrength}
    \mathbf{H}^s(\mathbf{r},t) = \nabla \times \iiint_{V_\varepsilon} \frac{\frac{\partial}{\partial t} \mathbf{J}_\varepsilon(\mathbf{r}',\tau)}{4 \pi R} \mathrm{d}V',
\end{equation}
from which we can derive the scattered electric field $\mathbf{E}^s$, i.e.
\begin{equation}
     \varepsilon_0\frac{\partial}{\partial t} \mathbf{E}^s (\mathbf{r},t)  =  - \frac{\partial}{\partial t}\mathbf{J}_\varepsilon(\mathbf{r},t) + \nabla \times \mathbf{H}^s.
\end{equation}
The total electric field is then the superposition of the incident and scattered electric fields, i.e.  $\mathbf{E}(\mathbf{r},t) = \mathbf{E}^i(\mathbf{r},t) + \mathbf{E}^s(\mathbf{r},t)$, which also determines the contrast current density as~\cite{Jin2010}
\begin{equation} \label{eq:ContrastCurrentDensity}
    \mathbf{J}_{\varepsilon}(\mathbf{r},t) = (\varepsilon_r(\mathbf{r})-1) \varepsilon_0 \mathbf{E}(\mathbf{r},t).
\end{equation}
By combining the above equations, removing the time derivative over the contrast current density and normalizing for the relative permittivity, we obtain the time domain contrast current density volume integral equation (TDJVIE)
\begin{equation} \label{eq:TDJVIE}
    \frac{\varepsilon_r(\mathbf{r})  -1}{\varepsilon_r(\mathbf{r})}   \varepsilon_0 \mathbf{E}^i(\mathbf{r},t)  =   \mathbf{J}_\varepsilon (\mathbf{r},t) -  \frac{\varepsilon_r(\mathbf{r})  -1}{\varepsilon_r(\mathbf{r})}\mathcal{S}(\mathbf{J}_\varepsilon)(\mathbf{r},t),
\end{equation}
with
\begin{equation}
    \mathcal{S}(\mathbf{F}) (\mathbf{r},t) =  \nabla \times \nabla \times \iiint_{V_{\varepsilon}} \frac{\mathbf{F}(\mathbf{r}',\tau)}{4 \pi R} \mathrm{d}V'.
\end{equation}
This definition of TDJVIE deviates from that in~\cite{VanDiepen2024}, but has the advantage that it does not require the time derivative of the incident electric field.

\subsection{Voxelization} \label{sc:Voxelization} 
Before discretizing the TDJVIE~\eqref{eq:TDJVIE} to find a numerical approximation to the contrast current density, we start by creating a piecewise-constant approximation of the relative permittivity of the scatterer $\varepsilon_r(\mathbf{r})$. To illustrate the step-by-step process, we apply the discretization to the dielectric sphere shown in Figure~\ref{fig:Sphere}, which has a constant $\varepsilon_r$ represented by the color red. First, we enclose the scatterer in a box, represented by the black dashed lines in Figure~\ref{fig:Sphere}, and divide the box evenly along each spatial Cartesian dimension, i.e. $U$-times in $\hat{\mathbf{x}}$-direction, $V$-times in $\hat{\mathbf{y}}$-direction, $W$-times in $\hat{\mathbf{z}}$-direction. This division results in a set of $M = U \times V \times W$ voxels all of dimension $\Delta x \times \Delta y \times \Delta z$, as shown in Figure~\ref{fig:Voxelizev3}. Owing to the regularity of the discretization, we differentiate between voxels using their index $m$ and integer index $[u,v,w]$ where
\begin{equation}
    m = (w-1)UV + (v-1)U + u
\end{equation}
with $u = 1,\ldots,U$, $v=1,\ldots,V$ and $w = 1,\ldots,W$ and therefore $m = 1,\ldots,M$. So, the $m$-th voxel has a corresponding integer index $[u,v,w]$ and occupies the volume $\mathcal{V}_{m}$, which is a beam of dimensions $\Delta x \times \Delta y \times \Delta z$ centered at the Cartesian coordinate $\mathbf{r}_{m} = \left((u-\frac{1}{2})\Delta x, (v-\frac{1}{2})\Delta y,(w-\frac{1}{2})\Delta z\right)$. The relative permittivity $\varepsilon_{m}$ throughout each $m$-th voxel is then equal to the relative permittivity of the non-discretized scattering setup at the location of the voxel center, i.e. $\varepsilon_{m} = \varepsilon_r(\mathbf{r}_{m})$. 
\begin{Figure}
    \centering
    \includegraphics[width = \columnwidth]{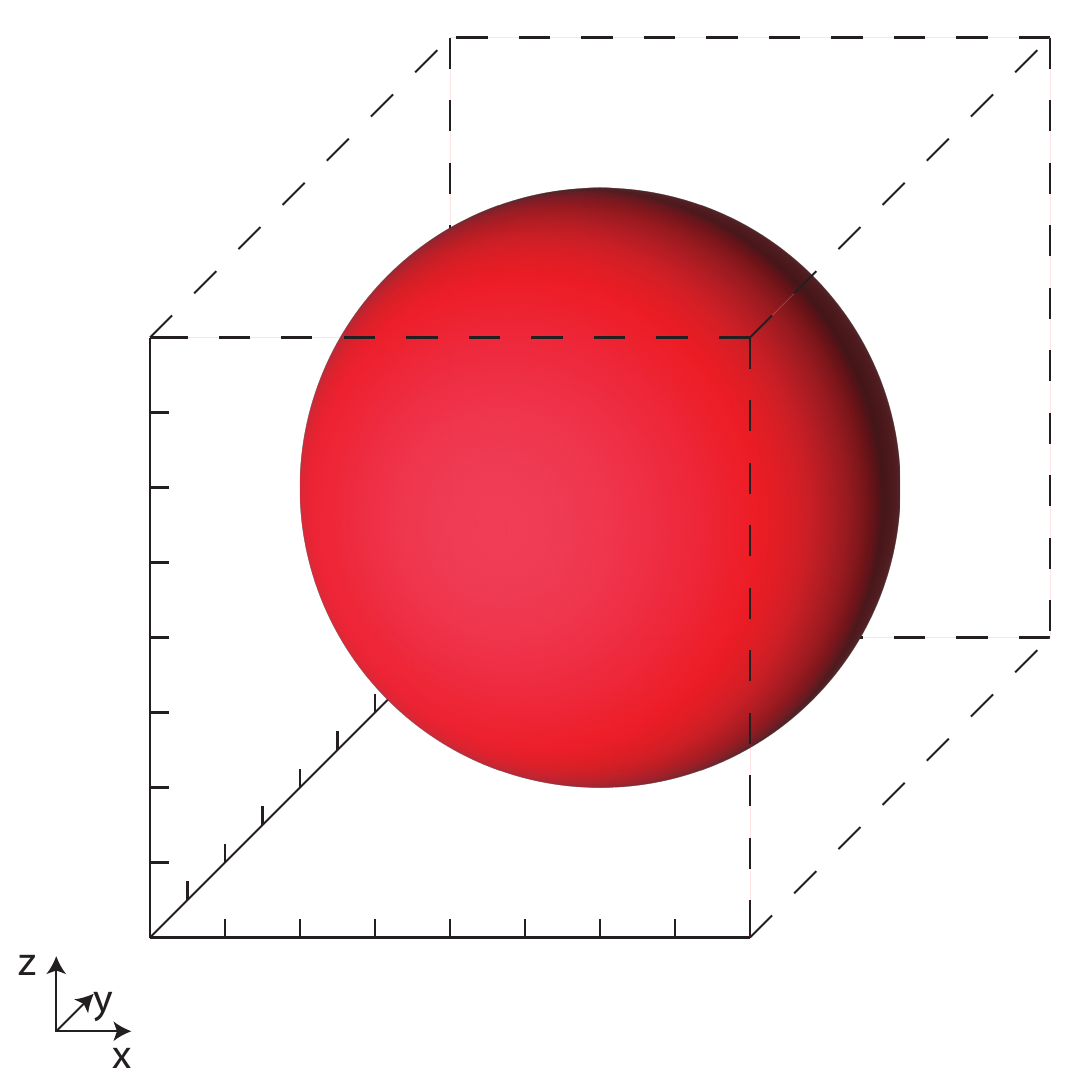}
    \captionof{figure}{The scattering setup consists of a sphere with a constant dielectric contrast value represented by the color red in a background medium. The sphere is enclosed by a box illustrated by the black dashed lines.}
    \label{fig:Sphere}
\end{Figure}
\begin{Figure}
    \centering
    \includegraphics[width = \columnwidth]{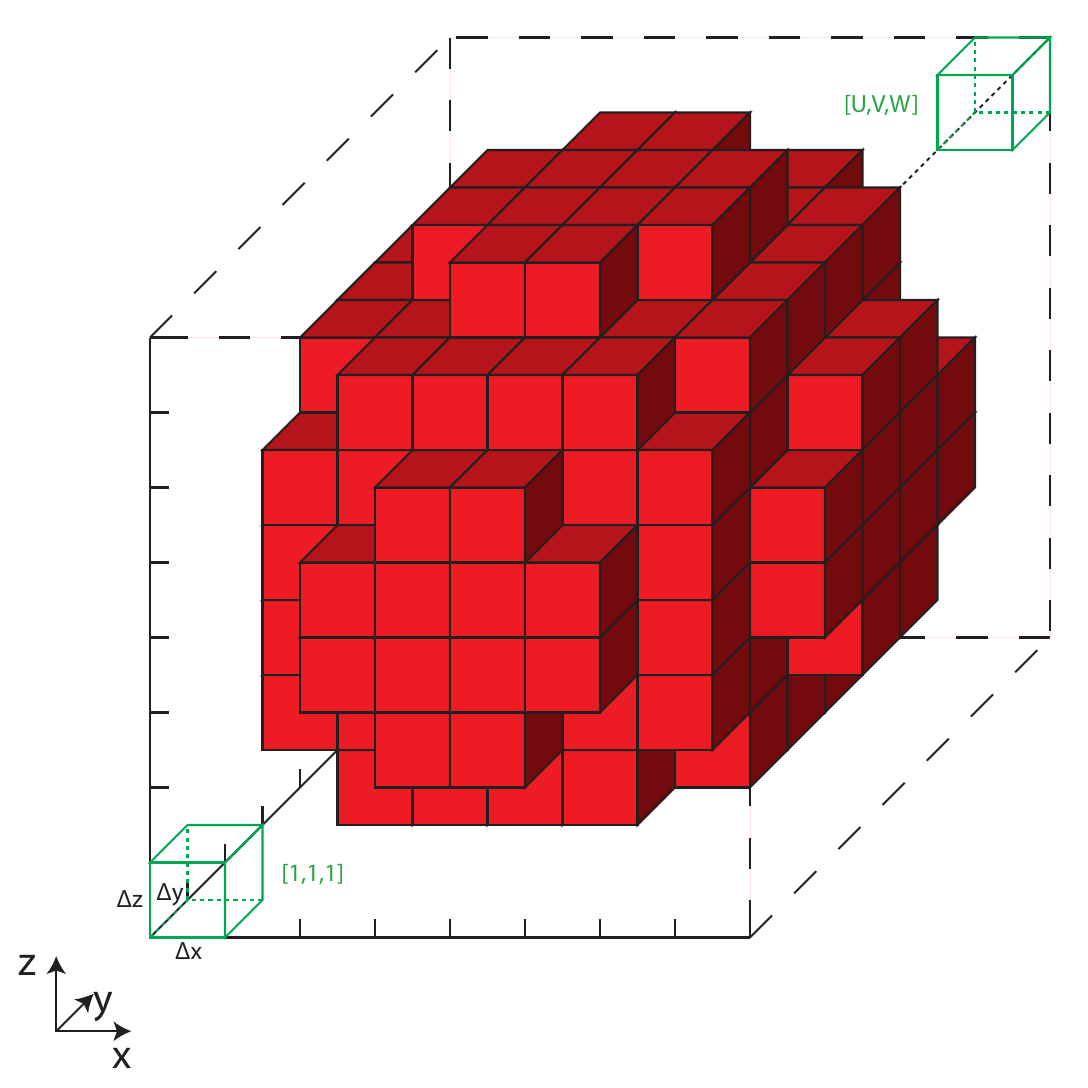}
    \captionof{figure}{The box around the scatter in Figure~\ref{fig:Sphere} is divided evenly along each Cartesian direction, i.e. $U$ times in the $\hat{\mathbf{x}}$-direction, $V$ times in the $\hat{\mathbf{y}}$-direction, $W$ times the in $\hat{\mathbf{z}}$-direction. So, the scattering setup is discretized using $M = U \times V \times W$ voxels, each of a dimension $\Delta x \times \Delta y \times \Delta z$. We have visualized the voxels in Green at $[1,1,1]$ and $[U,V,W]$. The dielectric contrast of each voxel is equal to the dielectric contrast of the scatterer setup shown in Figure~\ref{fig:Sphere} sampled at $\mathbf{r}_{\mathbf{u}}$, i.e. the center of each voxel.}
    \label{fig:Voxelizev3}
\end{Figure}

\subsection{Discretization} \label{sc:Discretization} 
As we now have a voxel representation of the scattering setup, we can move on to discretizing the TDJVIE~\eqref{eq:TDJVIE}. We employ the same discretization as introduced in~\cite{VanDiepen2024}, but now on a regularized voxel grid. We expand the dielectric contrast current density as
\begin{equation} \label{eq:J_basis}
    \mathbf{J}_\varepsilon(\mathbf{r},t) = \sum_{\alpha = x,y,z} \sum_{m' = 1}^{M} \sum_{n'=1}^N J^\alpha_{m',n'} \mathbf{f}^{\alpha}_{m'}(\mathbf{r}) T_{n'}(t)
\end{equation}
and introduce the testing operator on a vector field $\mathbf{g}$ as
\begin{equation}  
\label{eq:TestOperator}
    \mathcal{T}^\beta_{m,n}(\mathbf{g}) = \frac{1}{V} \int \delta (t-n\Delta t) \iiint \mathbf{f}^\beta_{m}(\mathbf{r}) \cdot \mathbf{g}(\mathbf{r},t) \mathrm{d}V \mathrm{d}t ,
\end{equation}
for $m = [1,1,1],\ldots,[U,V,W]$, $n = 1,\ldots,N$ and $\beta = x,y,z$. The variable $J_{m',n'}^\alpha$ in Equation~\eqref{eq:J_basis} is an expansion coefficient. The testing operator in Equation~\eqref{eq:TestOperator} works on a general three-dimensional vector field $\mathbf{g}(\mathbf{r},t)$ that depends both on space and time and contains the $\cdot$~operation, i.e. the scalar dot product between two three-dimensional vector functions. Further, $\mathrm{d}t$ is the infinitesimal time element and $\mathrm{d}V$ is the infinitesimal volume element over the observer coordinates. The functions $\mathbf{f}_{m'}^\alpha$~\eqref{eq:J_basis} and $\mathbf{f}_{m}^\beta$~\eqref{eq:TestOperator} are the spatial basis and test function, respectively. On each voxel $m'$ with $[u',v',w']$ in the voxelized box, we define three piece-wise constant basis functions, one for each spatial Cartesian direction, i.e. for $\alpha = x,y,z$ 
\begin{equation} \label{eq:SpatialBasis}
    \mathbf{f}^\alpha_{m'} (\mathbf{r}) = \begin{cases}
    \hat{\boldsymbol{\alpha}}, & \text{for} \, \mathbf{r} \in \mathcal{V}_{m'}, \\
    \mathbf{0}, & \text{for} \, \mathbf{r}  \notin \mathcal{V}_{m'}, \\
    \end{cases}
\end{equation}
where $\mathcal{V}_{m'}$ is the volume occupied by a voxel of dimension $\Delta x \times \Delta y \times \Delta z$ centered around $\mathbf{r}_{m'}$. In a similar way, we define $\mathbf{f}_{m}^\beta$ on the voxel centered around $\mathbf{r}_{m}$ occupying the volume $\mathcal{V}_{m}$, where $\mathbf{f}_{m}^\beta = \mathbf{f}_{m'}^\alpha$ if $\beta = \alpha$ and $m = m'$. The testing operator normalizes for the volume of the test voxel $V = \Delta x \Delta y \Delta z$. The functions $\delta(t-n \Delta t)$ and $T_{n'}(t)$ are the $n$-th temporal test and $n'$-th temporal basis function, respectively, with discrete time step size $\Delta t$. The stability of the resulting MOT-scheme depends on the choice of temporal basis and test functions~\cite{Sayed2015,VanDiepen2024}. The Dirac-delta test as $\delta(t-n \Delta t)$ in combination with the quadratic spline basis as $T_{n'}(t)$ results in a MOT-scheme where the stability does not depend on the dielectric contrast~\cite{VanDiepen2024}. The definition of the quadratic temporal spline basis function can be found in~\cite{VanDiepen2024}.

By substituting the contrast current density expansion~\eqref{eq:J_basis} in the TDJVIE~\eqref{eq:TDJVIE} and testing the resulting equation with the test operator in~\eqref{eq:TestOperator}, we end up with a matrix equation with interaction-matrix elements
\begin{equation} \label{eq:Z}
    Z^{\beta,\alpha}_{m,m',n,n'} =  \mathcal{T}^\beta_{m,n}\left( \mathbf{f}^{\alpha}_{m'}T_{n'} -\frac{\varepsilon_{m}-1}{\varepsilon_{m}}\mathcal{S}(\mathbf{f}^{\alpha}_{m'} T_{n'})\right),
\end{equation}
and excitation-vector elements
\begin{equation} \label{eq:E}
    E^\beta_{m',n'} = \mathcal{T}^\beta_{m,n} \left( \frac{\varepsilon_{m}-1}{\varepsilon_{m}}\varepsilon_0\mathbf{E}^i(\mathbf{r},t)\right),
\end{equation}
and unknown contrast-current-density vector elements $J_{m,n}$. The computation of the integrals in~\eqref{eq:Z} and~\eqref{eq:E} is discussed in~\cite{VanDiepen2024} and we adopt the same semi-analytic evaluation.

\subsection{MOT-scheme} \label{sc:MOTscheme}
In the MOT-scheme we exploit the discrete translation symmetry in time of the interaction matrix elements owing to the uniform expansion and sampling in time, i.e. the values of the matrix elements $Z^{\beta,\alpha}_{m,m',n,n'}$~\eqref{eq:Z} do not change if $n-n'$ does not change. The MOT interaction matrices $\mathbf{Z}_{n-n'}$ are created from computing the elements $Z^{\beta,\alpha}_{m,m',n-n'}$ for $n-n' = 0,\ldots,\ell$, $m = 1,\ldots,M$, $m'= 1,\ldots,M$, $\beta = x,y,z$ and $\alpha = x,y,z$. The total number of unique interaction matrices $ \ell$ is given by
\begin{equation} \label{eq:Ell}
    \ell = \left \lfloor \frac{\sqrt{(U \Delta x)^2+(V\Delta y)^2 + (W \Delta z)^2}}{c \Delta t} \right \rfloor + 2.
\end{equation} 
Together, these $\ell$ interaction matrices will form a banded lower block-triangular matrix equation of the form
\begin{equation} \label{eq:MatrixEquation}
    \begin{bmatrix}
    \mathbf{Z}_0 & \\
    \mathbf{Z}_1 & \mathbf{Z}_0 & & & \makebox(0,0){\text{\Huge0}} \\
    \vdots & \ddots & \ddots \\
    \mathbf{Z}_\ell & \cdots & \mathbf{Z}_1 &  \mathbf{Z}_0 \\
    \makebox(0,-10){\text{\huge0}} & \ddots & & \ddots & \ddots\\
    & & \mathbf{Z}_\ell & \cdots & \mathbf{Z}_1 &  \mathbf{Z}_0
    \end{bmatrix} 
    \begin{bmatrix}
    \mathbf{J}_1 \\
    \mathbf{J}_2 \\
    \vdots \\
    \mathbf{J}_\ell \\
    \vdots \\
    \mathbf{J}_n
    \end{bmatrix} = \begin{bmatrix}
    \mathbf{E}_1 \\
    \mathbf{E}_2 \\
    \vdots \\
    \mathbf{E}_\ell \\
    \vdots \\
    \mathbf{E}_n
    \end{bmatrix},
\end{equation}
where
\begin{equation} \label{eq:CurrentVector}
    \mathbf{J}_n = \left[J_{1,n}^x;J_{1,n}^y;J_{1,n}^z,\ldots;J_{M,n}^x;J_{M,n}^y;J_{M,n}^z\right],
\end{equation}
contains the expansion coefficients~\eqref{eq:J_basis} and
\begin{equation} \label{eq:ExcitationVector}
    \mathbf{E}_n = \left[E_{1,n}^x;E_{1,n}^y;E_{1,n}^z,\ldots;E_{M,n}^x;E_{M,n}^y;E_{M,n}^z\right],
\end{equation}
contains the excitation vector elements~\eqref{eq:E}. The semicolon $;$ that separate the elements in these and following expressions indicate that the elements form a column vector. As the matrix equation has a lower banded triangular matrix, forward substitution is used to find the solution at each time step $n$, i.e.
\begin{equation} \label{eq:MOTInversion}
    \mathbf{Z}_0 \mathbf{J}_n = \mathbf{E}_{n} - \mathbf{P}_{n-1},
\end{equation}
where
\begin{equation} \label{eq:PnMOT}
    \mathbf{P}_{n-1} = \sum_{n'=n-\ell}^{n-1} \mathbf{Z}_{n-n'} \mathbf{J}_{n'},
\end{equation}
with 
\begin{equation} 
    \mathbf{P}_n = \left[P_{1,n}^x;P_{1,n}^y;P_{1,n}^z;\ldots;P_{M,n}^x;P_{M,n}^y;P_{M,n}^z\right].
\end{equation}
The scheme in~\eqref{eq:MOTInversion} is the aforementioned MOT-scheme. The computational complexity of the MOT-scheme scales as $\mathcal{O}(M^2)$, due to the computation of $\mathbf{P}_{n-1}$~\eqref{eq:PnMOT}. Each interaction matrix $\mathbf{Z}_{n-n'}$ is sparse, but together these matrices will fill approximately a full matrix of $M^2$ elements. Thus, the MOT-scheme scales as $\mathcal{O}(M^2)$ per time step $n$.

\section{Spatial 
 FFT-acceleration} \label{sc:SpatialAcceleration}
The spatial FFT-acceleration of the MOT-scheme focuses on accelerating the matrix vector product $\mathbf{Z}_{n-n'} \mathbf{J}_{n'}$ in computing $\mathbf{P}_{n-1}$ in Equation~\eqref{eq:PnMOT} and is derived from FFT-acceleration techniques like CGFFT~\cite{Catedra1989,Zwamborn1991} and has been implemented for time domain surface integral equations~\cite{Yilmaz2001}. A uniform expansion and sampling in space results in a translation symmetry in space of the interaction matrix elements, i.e. the value of $Z^{\beta,\alpha}_{m,m',n,n'}$~\eqref{eq:Z} does not change if $[u-u',v-v',w-w']$ does not change. Now, by separating $\mathbf{P}_n$ and $\mathbf{J}_{n}$ along their respective Cartesian directions, i.e.
\begin{equation} 
    \mathbf{P}^\beta_n = \left[P_{1,n}^\beta;P_{2,n}^\beta,\ldots;P_{M,n}^\beta \right],
\end{equation}
for $\beta = x,y,z$, and,
\begin{equation} 
    \mathbf{J}^\alpha_{n} = \left[J_{1,n}^\alpha;J_{2,n}^\alpha,\ldots;J_{M,n}^\alpha \right],
\end{equation}
for $\alpha = x,y,z$, we replace the operation in~\eqref{eq:PnMOT} with
\begin{equation}~\label{eq:PnFFTMOT}
    \mathbf{P}^{\beta}_{n-1} = \sum_{n'=n-\ell}^{n-1} \mathbf{Z}^{\beta,\alpha}_{n-n'} \mathbf{J}^{\alpha}_{n'},
\end{equation}
where $\mathbf{Z}^{\beta,\alpha}_{n-n'}$ is now a three-level block-Toeplitz matrix, where each level corresponds to one of the spatial Cartesian directions. As discussed in Appendix~\ref{ap:ToeplitzFFTAcceleration}, the matrix vector product concerning a three-level block-Toeplitz matrix can be accelerated with three dimensional FFTs, where the associated FFT and point-wise multiplication operations scale as $\mathcal{O}(M \log M)$ and $\mathcal{O}(M)$, with $M$ the total number of voxels. As a result of the summation in~\eqref{eq:PnFFTMOT}, the matrix vector product is performed $\ell$ times per time step $n$. By clever reuse of the already computed FFTs and by performing the IFFT after the summation~\cite{Yilmaz2001}, Equation~\eqref{eq:PnFFTMOT} scales as $\mathcal{O}(M\log M) + \mathcal{O}(\ell M)$  per time step $n$. If $U \approx V \approx W$, then $\ell \sim \mathcal{O}(M^\frac{1}{3})$ and consequently Equation~\eqref{eq:PnFFTMOT} scales as $\mathcal{O}(M\log M) + \mathcal{O}(M^{\frac{4}{3}})$ per time step $n$, where the latter term will be dominant for large $M$. So, the scaling of spatial accelerated MOT over regular MOT is improved to $\mathcal{O}(M^{\frac{4}{3}})$ instead of $\mathcal{O}(M^2)$  per time step $n$. Still, further improvement in the computational complexity is required to go to larger sets of voxels.

\section{Temporal FFT-acceleration} \label{sc:TemporalAcceleration}
The computation of the values $\mathbf{P}_{n-1}$ every time step $n$ as in Equation~\eqref{eq:PnMOT} is equivalent to the following block-lower-triangular-Toeplitz matrix vector product,
\begin{equation} \label{eq:PnMatrixVector}
    \begin{bmatrix}
    \mathbf{P}_1 \\
    \vdots \\
    \mathbf{P}_\ell \\
    \vdots \\
    \mathbf{P}_n
    \end{bmatrix} = \begin{bmatrix}
    \mathbf{Z}_1 & & & \makebox(0,-10){\text{\Huge0}} \\
    \vdots & \ddots \\
    \mathbf{Z}_\ell & \cdots & \mathbf{Z}_1 &  \\
    \makebox(0,-10){\text{\huge0}}  & \ddots & & \ddots & \\
    & & \mathbf{Z}_\ell & \cdots & \mathbf{Z}_1 
    \end{bmatrix} 
    \begin{bmatrix}
    \mathbf{J}_1 \\
    \vdots \\
    \mathbf{J}_\ell \\
    \vdots \\
    \mathbf{J}_n
    \end{bmatrix}.
\end{equation}
We refer to this matrix as the MOT-matrix. The MOT-matrix is a Toeplitz matrix, owing to the uniform expansion and sampling in time, as discussed in Section~\ref{sc:MOTscheme}. Temporal FFT-acceleration of the MOT-scheme focuses on acceleration of this matrix vector product in the direction of the time stepping. This acceleration does not depend here on the internal matrix structure of the interaction matrices, so we first explain the temporal FFT-acceleration for scalar interaction matrices and subsequently combine it with the spatial FFT-acceleration discussed in Section~\ref{sc:SpatialAcceleration}. The temporal FFT-acceleration of the computation of Equation~\eqref{eq:PnMatrixVector} with scalar interaction matrices is similar to the one presented in~\cite{Hairer1985}, where non-linear Volterra convolution equations are considered. A difference is that the MOT-matrix in Equation~\eqref{eq:PnMatrixVector} is banded, whereas the one in~\cite{Hairer1985} is not.

\subsection{Toeplitz division} \label{sc:ToeplitzDivision}
The FFT-acceleration of a Toeplitz matrix vector product is well known and in Appendix~\ref{ap:ToeplitzFFTAcceleration} we further elaborate on this. However, its application to Equation~\eqref{eq:PnMatrixVector} is not straightforward, even though the MOT-matrix is a Toeplitz. To perform FFT-acceleration, all vector elements involved have to be known, however in the MOT-scheme, as a result of Equation~\eqref{eq:MOTInversion}, the value of $\mathbf{J}_n$ is not known before we have computed $\mathbf{P}_{n-1}$. Consequently, the FFT-acceleration cannot include elements above the diagonal in the MOT-matrix in Equation~\eqref{eq:PnMatrixVector}. Thus, to apply FFT-acceleration to MOT-schemes, we have to divide the MOT-matrix into smaller Toeplitz blocks that do not include elements above the diagonal. As explained in Appendix~\ref{ap:ToeplitzFFTAcceleration}, we achieve the highest FFT-acceleration by creating Toeplitz matrices with maximum dimension that are approximately square.
\begin{Figure}
    \centering
    \includegraphics[width = \columnwidth]{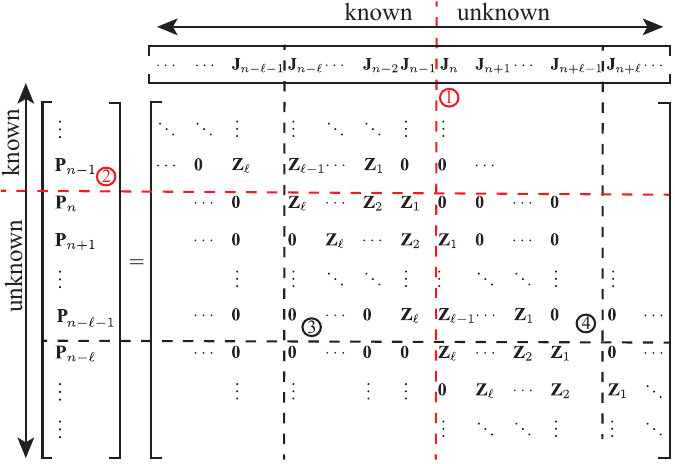}
    \captionof{figure}{The matrix elements shown here illustrate a part of the block-lower triangular matrix in Equation~\eqref{eq:PnMatrixVector}. We have computed the values of $\mathbf{J}_n$ up to the red-dashed line {\color{red} $\protect\numcircle{1}$}, thus we can perform any on the left of the line {\color{red} $\protect\numcircle{1}$}. Because we have computed up to $\mathbf{J}_n$, according to Equation~\eqref{eq:MOTInversion}, we have already solved all matrix vector products above the line {\color{red} $\protect\numcircle{2}$} to find $\mathbf{P}_{n-1}$. The optimally largest Toeplitz matrix $\mathbf{U}_0$ is the one enclosed by the horizontal and vertical red lines and the horizontal and vertical black dashed lines marked as $\protect\numcircle{3}$. The strictly lower triangular Toeplitz matrix $\mathbf{L}_0$ is the matrix enclosed by the horizontal and vertical red-dashed lines and the horizontal and vertical black-dashed lines marked as $\protect\numcircle{4}$}
    \label{fig:ToeplitzDivision}
\end{Figure}

We explain how to obtain the largest possible square Toeplitz matrix in Equation~\eqref{eq:PnMatrixVector} with the help of Figure~\ref{fig:ToeplitzDivision}, which represents a part of the Toeplitz matrix in Equation~\eqref{eq:PnMatrixVector}. In this figure, the red-dashed vertical line {\color{red} $\protect\numcircle{1}$} represents the $n$-th step in the MOT-scheme~\eqref{eq:MOTInversion}, i.e. everything left of this line concerns current densities that have already been computed, i.e. $\mathbf{J}_{n'}$ for $n' \leq n-1$, and everything on the right-hand side of this line concerns the current densities that are still to be computed and concern future time steps, i.e. $\mathbf{J}_{n'}$ for $n' \geq n$. At the $n$-th time step, we have already computed $\mathbf{P}_{n'}$ for $n' \leq n-1$, i.e. everything above the horizontal red-dashed line {\color{red} $\protect\numcircle{2}$}. So, we can form a Toeplitz matrix left of the vertical red-dashed {\color{red} $\protect\numcircle{1}$} and below the horizontal red-dashed line {\color{red} $\protect\numcircle{2}$}, for which the matrix vector product can be accelerated via FFTs owing to the fact that such a matrix only concerns current densities computed at earlier time steps. The optimally largest Toeplitz matrix to form is the one enclosed by the horizontal and vertical red lines and the horizontal and vertical black dashed lines marked as $\protect\numcircle{3}$. Any larger matrix is still Toeplitz but it will only include more rows and columns filled with exclusively zeros, which leads to larger FFT sizes without producing a more efficient matrix vector product. We define this matrix as
\begin{equation} \label{eq:ToeplitzU}
    \mathbf{U}_0 = \begin{bmatrix}
        \mathbf{Z}_\ell & \cdots &  \mathbf{Z}_1 \\
        & \ddots & \vdots \\
        & & \mathbf{Z}_\ell
    \end{bmatrix}.
\end{equation}
The product $\mathbf{U}_0 [\mathbf{J}_{n-\ell};\ldots;\mathbf{J}_{n-1}]$ is enough to compute $\mathbf{P}_n$, but for computing $\mathbf{P}_{n+1}$ up to $\mathbf{P}_{n+\ell-1}$ we also need the strictly lower triangular Toeplitz matrix to the right of $\mathbf{U}_0$, which is the matrix enclosed by the horizontal and vertical red-dashed lines and the horizontal and vertical black-dashed lines marked as $\protect\numcircle{4}$. We define this strictly lower triangular matrix as
\begin{equation} \label{eq:LowerTriangularToeplitzL}
    \mathbf{L}_0 = \begin{bmatrix}
        \mathbf{0} & \cdots & \cdots & \mathbf{0}\\
        \mathbf{Z}_1 &  \mathbf{0} & \cdots &\mathbf{0}\\
        \vdots & \ddots & \ddots & \vdots \\
        \mathbf{Z}_{\ell-1} & \cdots & \mathbf{Z}_1 & \mathbf{0}
    \end{bmatrix}.
\end{equation}
So, defining the largest possible Toeplitz matrix in Equation~\eqref{eq:PnMatrixVector} has resulted in a Toeplitz matrix $\mathbf{U}_0$ and strictly lower triangular Toeplitz matrix $\mathbf{L}_0$. The matrix shown in Figure~\ref{fig:ToeplitzDivision} can therefore also be written as
\begin{equation} \label{eq:ULDivision}
    \begin{bmatrix}
        \ddots & \ddots \\
            & \mathbf{U}_0 & \mathbf{L}_0 \\
            & & \mathbf{U}_0 & \mathbf{L}_0 \\
            & & & \mathbf{U}_0 & \mathbf{L}_0 \\
            & & & & \ddots & \ddots 
    \end{bmatrix}.
\end{equation}
We can now formally write the computation of vector $\mathbf{P}_n$ at time steps $(n'-1)\ell+1$ to $n'\ell$ as the matrix vector product involving two matrices, i.e.
\begin{equation} \label{eq:ULMatrixVectorProduct}
    \begin{bmatrix}
        \mathbf{P}_{(n'-1)\ell+1} \\
        \vdots \\
        \mathbf{P}_{n'\ell} \\
    \end{bmatrix} = 
    \mathbf{U}_0
    \begin{bmatrix}
        \mathbf{J}_{(n'-2)\ell+1} \\
        \vdots \\
        \mathbf{J}_{(n'-1)\ell} \\
    \end{bmatrix}
    +
    \mathbf{L}_0
    \begin{bmatrix}
        \mathbf{J}_{(n'-1)\ell+1} \\
        \vdots \\
        \mathbf{J}_{n'\ell} \\
    \end{bmatrix},
\end{equation}
where the matrix vector product with $\mathbf{U}_0$ can be computed via FFT-acceleration as it only involves $\mathbf{J}_n$ at previous time steps that have already been computed. The matrix vector product involving $\mathbf{L}_0$ cannot directly be computed via FFT-acceleration, as it involves the solution at future time steps. However, as $\mathbf{L}_0$ is a strictly lower triangular Toeplitz matrix, just like the original MOT-matrix in Equation~\eqref{eq:PnMatrixVector}, we can apply the same steps as for Equation~\eqref{eq:PnMatrixVector}. This is visualized in Figure~\ref{fig:LDivision}, where we divide $\mathbf{L}_0$ into the elements above the diagonal that we cannot include in the Toeplitz, marked by the upper-right red triangle, the optimally largest Toeplitz matrix $\mathbf{U}_1$ for which we can apply FFT-acceleration and two smaller strictly lower triangular matrices $\mathbf{L}_1$. We can repeat this process recursively for all $\mathbf{L}_k$ until some $k = K$ where $\mathbf{U}_K = \mathbf{Z}_1$, so further division is not possible anymore. The approximate dimensions of each square Toeplitz matrix $\mathbf{U}_k$ half with each iteration as shown in Figure~\ref{fig:LDivision}. However, their actual dimension should be an integer, but this will be addressed in Section~\ref{sc:SpatialTemporalAcceleration}. Overall, as the dimension of $\mathbf{U}_k$ halves with each iteration, the number of unique Toeplitz matrices $\mathbf{U}_k$ scales as $\mathcal{O}(\log \ell)$.
\begin{Figure}
    \centering
    \includegraphics[width = \columnwidth]{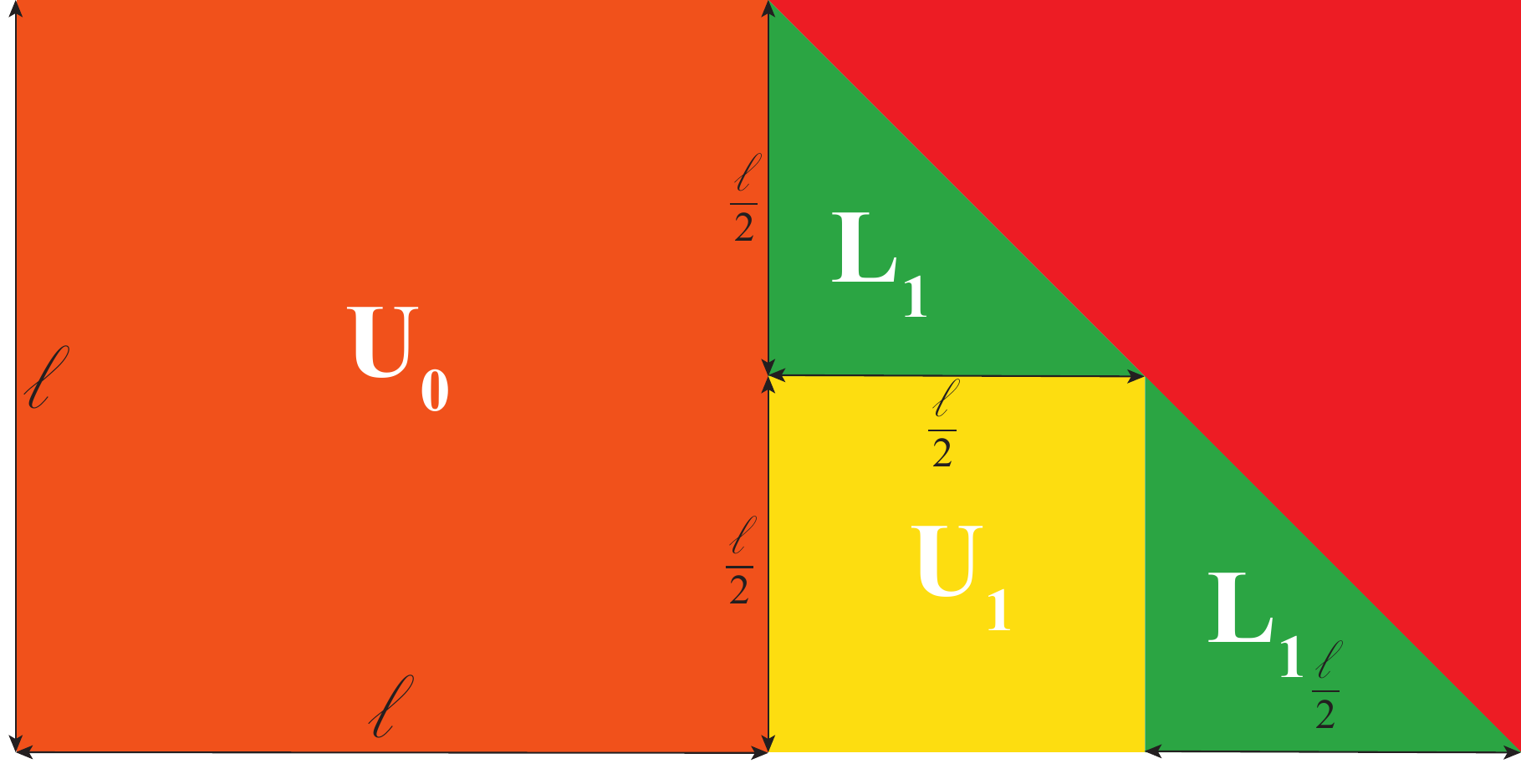}
    \captionof{figure}{The division of the strictly lower triangular Toeplitz matrix $\mathbf{L}_0$~\eqref{eq:LowerTriangularToeplitzL} into strictly lower triangular Toeplitz matrices $\mathbf{L}_1$ and Toeplitz matrix $\mathbf{U}_1$ to which we can apply FFT-acceleration as it does not include elements above the diagonal in Equation~\eqref{eq:PnMatrixVector} here marked by the red triangle. We indicate the approximate sizes of the matrices.}
    \label{fig:LDivision}
\end{Figure}

\subsection{Complexity} \label{sc:Complexity}
In Section~\ref{sc:ToeplitzDivision} we explained how we divide the banded Toeplitz matrix in Equation~\eqref{eq:PnMatrixVector} into smaller blocks of Toeplitz matrices $\mathbf{U}_k$ to which we can apply FFT-acceleration. This division divides the whole MOT-matrix into portions that cover $\ell$ time steps, see Equation~\eqref{eq:ULMatrixVectorProduct}. Thus the total complexity of temporal accelerated MOT-scheme per time step depends on what happens in these $\ell$ time steps and dividing that by $\ell$.

The dimension of $\mathbf{U}_k$ halves with each level $k$, i.e. $\mathbf{U}_0$ has dimension $\ell \times \ell$ and consequently $\mathbf{U}_k$ has dimension $\frac{\ell}{2^k} \times \frac{\ell}{2^k}$ as illustrated in Figure~\ref{fig:LDivision}. Because the dimension halves with each level $k$, the occurrence of $\mathbf{U}_k$ in $\ell$ times steps increases, i.e. $\mathbf{U}_k$ for $k \geq 1$ occurs $2^{k-1}$ times. The complexity of the individual matrix vector product with $\mathbf{U}_k$ scales as $\mathcal{O}(\frac{\ell}{2^k} \log \frac{\ell}{2^k})$. Thus, the complexity of the combined matrix vector products concerning $\mathbf{U}_k$ for $k \geq 1$ per time step scales as
\begin{equation} \label{eq:ComplexityHierarchicalDivisionInTime}
    \sum_{k=1}^K \frac{1}{\ell} 2^{k-1} \mathcal{O}\left(\frac{\ell}{2^k}  \log \frac{\ell}{2^k}\right) = \mathcal{O}\left(K \log \frac{\ell}{2^{\frac{K+1}{2}}} \right).
\end{equation}
As explained at end of Section~\ref{sc:ToeplitzDivision}, $K$ scales as $\mathcal{O}( \log \ell)$. Consequently, the average complexity per time step of the MOT-scheme with temporal FFT-acceleration scales as $\mathcal{O}( \log^2 \ell)$, which is in agreement with literature~\cite{Hairer1985}.

\section{Spatial-temporal FFT-acceleration} \label{sc:SpatialTemporalAcceleration}
The spatial FFT-acceleration of the MOT-scheme discussed in Section~\ref{sc:SpatialAcceleration} and the temporal FFT-acceleration of the MOT-scheme discussed in Section~\ref{sc:TemporalAcceleration} can now be combined. The uniform expansion and sampling in space and time results in a translation symmetry in space and time of the interaction matrix elements, i.e. the value of $Z^{\beta,\alpha}_{m,m',n,n'}$~\eqref{eq:Z} does not change if $[u-u',v-v',w-w']$ and $[n-n']$ do not change. Therefore, the $K \sim \log \ell$ unique block-Toeplitz matrices $\mathbf{U}_k$ introduced in Section~\ref{sc:ToeplitzDivision} are each a four-level block Toeplitz matrix of dimension $\mathcal{O}(\frac{\ell}{2^k} M)$. As explained at the end of Appendix~\ref{ap:ToeplitzFFTAcceleration}, the matrix vector product concerning $\mathbf{U}_k$ can then be accelerated with four dimensional FFTs whose complexity scales as $\mathcal{O}(M \frac{\ell}{2^k} \log (M \frac{\ell}{2^k}))$. As discussed in Section~\ref{sc:Complexity}, the larger matrices $\mathbf{U}_k$ have a lower occurrence in time and therefore the average complexity per time step of the spatial-temporal accelerated MOT-scheme scales as $\mathcal{O}(M \log \ell \log (M \ell))$. The values of $\ell$ and $M$ are coupled for time domain integral equations due to the propagating Green function~\eqref{eq:GreenFunction}. Consequently, $\ell$ scales as $\mathcal{O}(M^\frac{1}{d})$, see Equation~\eqref{eq:Ell}, where $d$ represents the physical dimension of the scatterer, which can be one, two or three dimensional. Independent of the dimension, $\log \ell \sim \log M$ and the complexity of spatial-temporal accelerated MOT-scheme is then rewritten to $\mathcal{O}(M \log^2 M)$. This is the same as the complexity of the spatial-temporal FFT-acceleration presented in~\cite{Yilmaz2002a}, which employs the Toeplitz division as discussed in Section~\ref{sc:ToeplitzDivision} for the time domain surface integral equations. Starting at $\mathbf{U}_K = \mathbf{Z}_1$ and defining the rest from there, alleviates the issue of  of $\frac{\ell}{2^k}$ not being an integer.  

The Toeplitz divison as discussed in Section~\ref{sc:ToeplitzDivision} hierarchically divides time to obtain the Toeplitz matrices $\mathbf{U}_k$. In the case of time domain integral equations space and time are coupled due to the Green function~\eqref{eq:GreenFunction}. Consequently, the matrix $\mathbf{U}_k$ is limited to interaction between basis and test elements separated no more than a predefined radial distance $R \leq (2^{K-k}+2) c\Delta t$. Therefore, the hierarchical division of time to obtain $\mathbf{U}_k$ divides the mesh into $K$ levels, where $\mathbf{U}_K$ only includes interactions close by, but $\mathbf{U}_{K-1}$ includes interactions at double that distance, and $\mathbf{U}_k$ at $2^{k-k}$-times that distance. Thus, a hierarchical division of time leads to a hierarchical division of space. The work in~\cite{Yilmaz2002} shows that it also works the other way around, i.e. a hierarchical division of space leads to a hierarchical division of time from where one can define Toeplitz matrices similar to $\mathbf{U}_k$. They refer to this technique as HIL-FFT and apply it to time domain surface integral equations to also achieve a $\mathcal{O}(M \log^2 M)$-scaling.

\subsection{3D spatial-temporal FFT-acceleration} \label{sc:3D}
We extend the HIL-FFT~\cite{Yilmaz2002}, i.e. spatial-temporal FFT-acceleration via a hierarchical division of space, to 3D. We define the four-level Toeplitz matrices $\mathbf{U}_k$ by dividing the interaction matrices into $K \sim \log \ell$ levels. Inspired by~\cite{Yilmaz2002a}, we start at the definition of $\mathbf{U}_K$, unlike~\cite{Yilmaz2002} which starts by defining $\mathbf{U}_0$. The four-level Toeplitz matrix $\mathbf{U}_K$ contains all interaction matrix elements in $\mathcal{M}_K$, where $\mathcal{M}_K$ contains all basis and test voxel pairs for which $-U_K \leq u-u' \leq U_K$, $-V_K \leq v-v' \leq V_K$ and $-W_K \leq w-w' \leq W_K$ holds. As explained in Section~\ref{sc:Complexity}, $\mathbf{U}_K$ will have the highest occurrence, thus its dimension should remain small, i.e. minimize $U_L$, $V_L$, and $W_L$ such that $\mathbf{Z}_1$ is included in $\mathbf{U}_K$ as it did in Section~\ref{sc:ToeplitzDivision}. The next four-level Toeplitz matrix $\mathbf{U}_{K-1}$ contains all interaction matrix elements in $\mathcal{M}_{K-1}$ but not in $\mathcal{M}_{K}$, i.e. $\mathcal{M}_{K-1} \backslash \mathcal{M}_{K}$, where $\mathcal{M}_{K-1}$ contains all basis and voxel pairs for which $-U_{K-1} \leq u-u' \leq U_{K-1}$, $-V_{K-1} \leq v-v' \leq V_{K-1}$ and $-W_{K-1} \leq w-w' \leq W_{K-1}$ holds. We repeat this process, where a level $k$ concerns the four-level Toeplitz matrix $\mathbf{U}_k$ which contains all interaction matrix elements in $\mathcal{M}_{k}$ but not in $\mathcal{M}_{k+1}$, i.e. $\mathcal{M}_{k} \backslash \mathcal{M}_{k+1}$, where $\mathcal{M}_{k}$ contains all basis and voxel pairs for which $-U_{k} \leq u-u' \leq U_{k}$, $-V_{k} \leq v-v' \leq V_{k}$ and $-W_{k} \leq w-w' \leq W_{k}$ holds. To obtain $K \sim \log \ell$ levels, we define the relation between $\mathcal{M}_k$ and $\mathcal{M}_{k+1}$ as $U_k = 2U_{k+1}$, $V_k = 2V_{k+1}$ and $W_k = 2W_{k+1}$ with a limit to these values $U_0 = U-1$, $V_0 = V-1$ and $W_0 = W-1$, i.e. the original size of the problem. This is similar to the doubling of the hierarchical division of time with a maximum dimension of $\ell$ as explained in Section~\ref{sc:ToeplitzDivision}. Thus, the definition of $\mathcal{M}_k$ for $k = 0,\ldots,K$, splits the interaction matrices into $K+1$ sets as illustrated in Figure~\ref{fig:SpatialHierarchicalDivision}.
\begin{Figure}
    \centering
    \includegraphics[width=0.8\columnwidth]{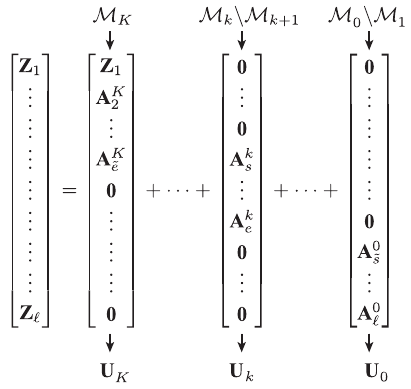}
    \captionof{figure}{The definition of the sets $\mathcal{M}_k$ is used to divided the interaction matrices $\mathbf{Z}_{1}$ through $\mathbf{Z}_\ell$ into $K$ sets. The $k$-th set is used to construct the four-level Toeplitz matrix $\mathbf{U}_k$.}
    \label{fig:SpatialHierarchicalDivision}
\end{Figure}
The minimization of $U_K$, $V_K$ and $W_K$ results in $\mathbf{A}^K_1 = \mathbf{Z}_1$. The leading zeros-matrices included in $\mathbf{U}_k$ in Figure~\ref{fig:SpatialHierarchicalDivision} before time step $s$ are a result of the travel time of the wave front from a basis voxel to the test voxels in $\mathcal{M}_k \backslash \mathcal{M}_{k+1}$ and the trailing zeros-matrices included in $\mathbf{U}_k$ in Figure~\ref{fig:SpatialHierarchicalDivision} after time step $e$ are a result of the back of that same wave leaving all test voxels in $\mathcal{M}_k \backslash \mathcal{M}_{k+1}$. Substituting $[\mathbf{A}^{k}_s;\ldots;\mathbf{A}^{k}_e]$ in Equation~\eqref{eq:PnMatrixVector}, we repeat the Toeplitz division as in Section~\ref{sc:ToeplitzDivision} and define the four-level Toeplitz matrix
\begin{equation}
    \mathbf{U}_k = \begin{bmatrix}
        \mathbf{A}^k_s & \mathbf{0} & \cdots & \mathbf{0}\\ 
        \vdots & \ddots & \ddots & \vdots  \\
        \vdots &  & \ddots & \mathbf{0} \\
        \mathbf{A}^k_{e}  & \cdots & \cdots & \mathbf{A}^k_s \\
        \mathbf{0} & \ddots & & \vdots \\
        \vdots & \ddots & \ddots & \vdots \\        
        \mathbf{0} & \cdots & \mathbf{0} & \mathbf{A}^k_{e}
    \end{bmatrix}.
\end{equation}
where the dimension is $U \cdot V \cdot W \cdot e \times U \cdot V \cdot W \cdot s$ and it occurs $\frac{\ell}{s}$-times in $\ell$ time steps.

To ease the complexity analysis, we consider the case $U \approx V \approx W \approx \sqrt[3]{M}$ and $\Delta x \approx \Delta y \approx \Delta z$, then $\ell \sim \sqrt[3]{M}$~\eqref{eq:Ell}.
Subsequently, both $e$ and $s$ are proportional to the dimensions of $\mathcal{M}_k$ and half for higher levels of $k$, i.e. $e \sim \mathcal{O}(\frac{\ell}{2^k})$ and $s \sim \mathcal{O}(\frac{\ell}{2^k})$. Consequently, the dimension of $\mathbf{U}_k$ scales as $\mathcal{O}(\frac{M \ell}{2^{k}})$ and it occurs $\mathcal{O}(2^k)$-times in $\ell$ time steps.
The complexity of the individual matrix vector product with $\mathbf{U}_k$ scales as $\mathcal{O}(\frac{M \ell}{2^{k}} \log \frac{M \ell}{2^{k}})$ as explained in Appendix~\ref{ap:ToeplitzFFTAcceleration}. Thus, the average complexity of the combined matrix vector products concerning $\mathbf{U}_k$ per time step scales as
\begin{equation}~\label{eq:ComplexityHierarchicalDivisionInSpace}
    \sum_{k=0}^{K} \frac{\mathcal{O} \!\left(2^k\right)}{\ell} \mathcal{O}\!\left(\frac{M \ell}{2^{k}} \log\frac{M \ell}{2^{k}} \right)  = \mathcal{O}\! \left(K M \log \frac{M \ell}{2^{\frac{K+1}{2}}} \right),
\end{equation}
where $M$ is the number of voxels used in the discretization of the MOT-JVIE. The complexity scaling of the hierarchical division in space~\eqref{eq:ComplexityHierarchicalDivisionInSpace} is thus similar to the complexity scaling of the hierarchical division of time~\eqref{eq:ComplexityHierarchicalDivisionInTime}, but with block interaction matrices. As $K \sim \log \ell$ and $\ell \sim \sqrt[3]{M}$, this reduces to the expected $\mathcal{O}(M \log^2 M)$. Starting at $\mathbf{U}_K$ instead of $\mathbf{U}_0$ as in~\cite{Yilmaz2002} has therefore not altered the complexity scaling. However, it prevents the need for $U$, $V$ and $W$ to be a multiple of two as we can easily truncate $\mathcal{M}_k$ to the required size. We also recommend a different implementation of the construction of $\mathbf{U}_k$. We compute the elements of $Z^{\beta,\alpha}_{m,m',n,n'}$~\eqref{eq:Z} for $-U+1 \leq u-u' \leq U-1$, $-V+1 \leq v-v' \leq V-1$ and $-W+1 \leq v-v' \leq W-1$ and $n-n' = 0,\ldots,\ell$, which we store in a four dimensional array from where we derive $\mathbf{U}_k$. This is significantly simpler to implement than the intensive bookkeeping proposed by~\cite{Yilmaz2002}.

\subsection{Numerical results} \label{sc:AccelerationNumericalResults}
To demonstrate the MOT-JVIE spatial-temporal FFT-acceleration discussed in Section~\ref{sc:3D}, we will compute the contrast current density inside a $0.2^3 ~\mathrm{m}^3$ cubic $\varepsilon_r = 12$ scatterer centered at $\mathbf{r} = (0.1,0.1,0.1)$, induced by a Gaussian $\hat{\mathbf{x}}$-polarized plane wave travelling in the negative $\hat{\mathbf{z}}$-direction defined as
\begin{equation} \label{eq:GaussPlaneWave}
    \mathbf{E}^i (\mathbf{r},t) = \frac{4 E_0 }{\sigma\sqrt{\pi}} \hat{\mathbf{p}} \exp{\left(-\left(\frac{4}{\sigma}((t-t_0) - \mathbf{r}\cdot \hat{\mathbf{k}}))\right)^2 \right)},
\end{equation}
with polarization $\hat{\mathbf{p}} = \hat{\mathbf{x}}$, propagation direction $\hat{\mathbf{k}} = -\hat{\mathbf{z}}$ and $E_0$ is the amplitude scaling, set to $E_0 = 1~\mathrm{V}/\mathrm{m}$.  The unit $\mathrm{lm}$ is known as lightmeter, i.e. the time it takes for the wave to travel a distance of $1~\mathrm{m}$ and is used in Equation~\eqref{eq:GaussPlaneWave} for the pulse width $\sigma$, set to $\sigma = 2~\mathrm{lm}$, and the separation time at time $t = 0$ between the Gaussian pulse center and the coordinate system origin $t_0$, set to $t_0 = 3.42~\mathrm{lm}$. To test the acceleration, we have to increase the number of voxels in the discretization $M$. As explained in Section~\ref{sc:Voxelization}, the voxels are defined by enclosing the scatterer by a box divided evenly along each Cartesian direction, where $U = V = W = 0.2/\sqrt[3]{M}$. Consequently, the dimensions of a voxel are equal, i.e. $\Delta x = \Delta y = \Delta z = 0.2/\sqrt[3]{M}$. We set $\Delta t = \Delta x /c_0$ to maintain: sparsity in the interaction matrix $\mathbf{Z}_0$; accuracy in the numerical evaluation of the volume test integral in Equation~\eqref{eq:TestOperator}~\cite{VanDiepen2024}; and  have $\ell \sim M^{1/3}$.

We have implemented the MOT-JVIE with MATLAB R2018b version and ran that on two Intel(R) Xeon(R) Gold 6148 CPU's @ 2.40GHz for $\sqrt[3]{M} \in \{20,36,52,64,100\}$. The average computation time to compute $\mathbf{P}_n$ per time step as a function of $M$ is shown in Figure~\ref{fig:Cube_epsr12_acceleration}.
\begin{Figure}
    \centering
    \includegraphics[width=\columnwidth]{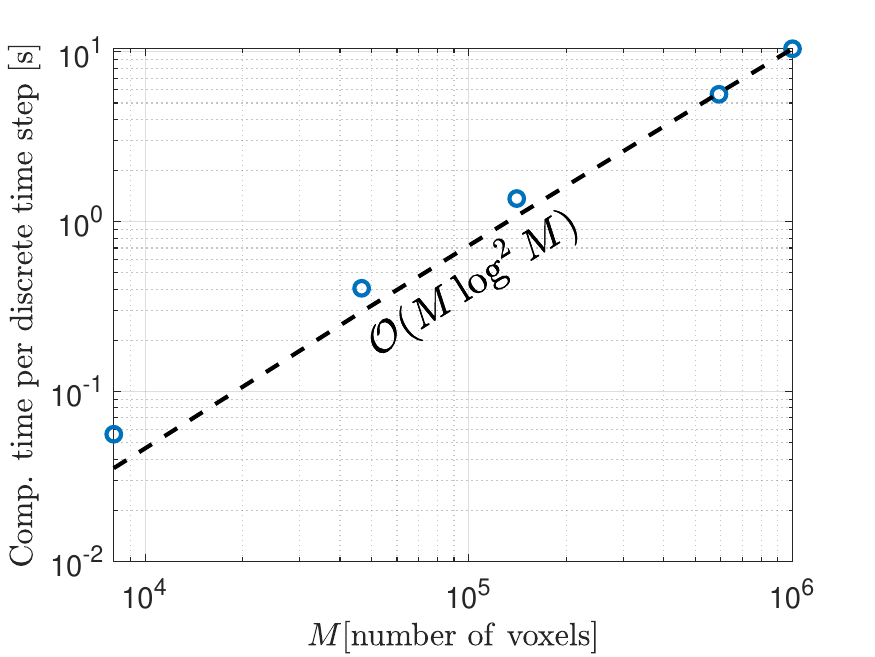}
    \captionof{figure}{The average computation time of $\mathbf{P}_n$ per time step in the spatial-temporal accelerated MOT-JVIE implemented as explained in Section~\ref{sc:3D} with MATLAB R2018b running on two Intel(R) Xeon(R) Gold 6148 CPU's @ 2.40GHz.}
    \label{fig:Cube_epsr12_acceleration}
\end{Figure}
The black-dashed line indicates the theoretical $\mathcal{O}(M \log^2 M)$-scaling, visually confirming the $\mathcal{O}(M \log^2 M)$-scaling in the spatial-temporal FFT-acceleration of the MOT-JVIE as explained in Section~\ref{sc:3D}.

To study the solution convergence for a higher number of voxels, we have sampled the contrast current density at $\mathbf{r} = (x,y,z)$ with $x \in \{0.025,\! 0.075,\!0.125,\!0.175\}$, $y \in \{0.025,\!0.075,\!0.125,\!0.175\}$ and $z \in \{0.025,\!0.075,\!0.125,\!0.175\}$  in this numerical experiment, i.e. $64$ locations in total. These sample locations remain at the center of the voxels, which is important for a convergence study~\cite{VanDiepen2024}. The $\hat{\mathbf{x}}$-component of the solution, $J^x(\mathbf{r},t)$, at $\mathbf{r} = (0.025,0.075,0.025)$ for the different $\sqrt[3]{M}$ are shown in Figure~\ref{fig:Cube_epsr12_InstabilityComparison}.
\begin{Figure}
    \centering
    \includegraphics[width=\columnwidth]{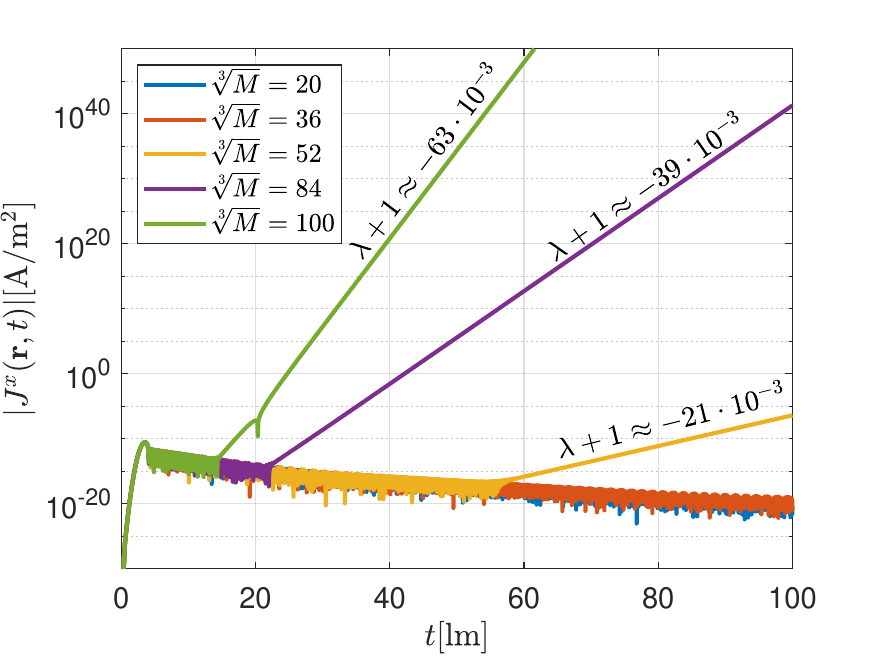}
    \captionof{figure}{The $\hat{\mathbf{x}}$-component of the MOT-JVIE solution, $J^x(\mathbf{r},t)$, as a function of time sampled at $\mathbf{r} = (0.025,0.075,0.025)$ in the $0.2^3~\mathrm{m}^3$ dielectric cube with $\varepsilon_r = 12$. The discretization settings are $\Delta x = \Delta y = \Delta z = c_0 \Delta t = 0.2/\sqrt[3]{M}~\mathrm{m}$. The companion matrix eigenvalues $\lambda$ associated to the nonphysical unstable solution are added to the respective lines.}
    \label{fig:Cube_epsr12_InstabilityComparison}
\end{Figure}
We observe a nonphysical unstable per-time-step alternating solution at the highest frequency $f = 1/(2\Delta t)$, i.e. the solution has a sign flip between discrete time steps and there is an exponential monotonic increase in magnitude, for $\sqrt[3]{M} \geq 52$. This type of instability is associated with the companion matrix eigenvalues, $\lambda$, on the negative real axis outside the unit circle~\cite{Sayed2015,VanDiepen2024,VanDiepen2024a}. We have estimated the eigenvalues from the solutions and added their values to the respective lines. We observe that the eigenvalue moves further away from the unit circle with an increase in $\sqrt[3]{M}$. We will address this instability in the next section.

\section{Stabilization} \label{sc:Stabilization}
The spatial-temporal FFT-acceleration of the MOT-JVIE presented in Section~\ref{sc:3D} enables simulations with a number of voxels above $M = 20^3$ within reasonable computation times, see Figure~\ref{fig:Cube_epsr12_acceleration}. The numerical experiments in Section~\ref{sc:AccelerationNumericalResults} illustrate that the MOT-JVIE suffers from an instability when the number of voxels increases.
The pertaining nonphysical unbounded solution corresponds to the eigenvalues of the companion matrix, $\lambda$, close to $-1$, but just outside the unit circle. We conjecture that the finite precision in the calculation of the interaction matrix elements accumulates to instability as the number of voxels in the simulation increases. To analyse this behavior, we employ the positive definite stability analysis (PDSA) presented in~\cite{VanDiepen2024a}.

The PDSA is a stability analysis technique derived from the companion matrix stability analysis. The PDSA guarantees that all eigenvalues of the companion matrix that lie on the negative real axis are within the unit circle if the matrices 
\begin{equation} \label{eq:Dn}
    \mathbf{D}_n = \binom{\ell}{\ell-n} \mathbf{Z}_0 + \ldots + (-1)^n \binom{\ell-n}{\ell-n} \mathbf{Z}_{n}
\end{equation}
for $n = 0,\ldots,\ell$ are all positive definite. Computing $\mathbf{D}_n$ in the numerical experiment of Section~\ref{sc:AccelerationNumericalResults} for smaller $\sqrt[3]{M}$ 
shows that $\mathbf{D}_\ell$ is not positive definite from $\sqrt[3]{M} \geq 14$. Losing the positive definiteness of $\mathbf{D}_\ell$ means that some of the negative real eigenvalues are potentially outside the unit circle and this is in line with the observations in Section~\ref{sc:AccelerationNumericalResults}. A possible reason for this is the accumulation of finite precision effects in the calculation of the interaction matrix elements. The error due to finite precision in an interaction matrix, represented as a matrix $\mathbf{N}$, leads to an offset in the matrix $\mathbf{D}_\ell$, i.e.
\begin{equation} \label{eq:DellFinitePrecision}
    \hat{\mathbf{D}}_\ell = \mathbf{D}_\ell + \mathbf{N}.
\end{equation}
These matrices are all symmetric because the operator in Equation~\eqref{eq:Z} is symmetric and the actual matrices are computed by exploiting this symmetry. The eigenvalues of symmetric matrices are real~\cite{VanDiepen2024a} and we can define the eigenvalue range as $\lambda(\mathbf{D}_\ell) \in [\lambda^-,\lambda^+]$ and $\lambda(\mathbf{N}) \in [\mu^-,\mu^+]$. The eigenvalues of the sum of two symmetric matrices are then bounded by sum of their ranges~\cite{Horn2013}, i.e. $\lambda(\hat{\mathbf{D}}_\ell) \in [\lambda^- + \mu^-,\lambda^+ + \mu^+]$. We can find a lower bound for $\mu^-$ from unifying the disks in the Gershgorin theorem~\cite{Horn2013}, which leads to $\mu^- \ge -3M\epsilon$, where $\epsilon$ is the largest absolute finite precision error in $\mathbf{N}$ and $3M$ is the number of unknowns. Consequently, if $3M <  \lambda^-/\epsilon$, the PDSA guarantees that the MOT-JVIE is stable in the presence of finite-precision effects in the elements of interaction matrices. The lower bound $\lambda^-$ is governed by the physics and the choice for the discretization. Thus, the number of voxels in a stable MOT-JVIE discretization is limited by the limited accuracy in the interaction matrix elements due to finite precision in the numerical calculation of the underlying integrals.

To estimate the value of $\epsilon$, we repeat the numerical experiment for $\sqrt[3]{M}=20$ in Section~\ref{sc:AccelerationNumericalResults}, but we make some alterations. First, we lower the permittivity of the cube to $\varepsilon_r = 2$. The error in the interaction matrix values $\epsilon$ is independent of the permittivity and lowering the permittivity will remove the resonances in our solution and lets us focus the nonphysical unstable per-time-step alternating solution. Second, we deliberately introduce an additional error $\epsilon_{\text{trunc}}$ by truncating the interaction matrix values accordingly. We increase $\epsilon_{\text{trunc}}$ in the truncated MOT-JVIE until we observe a difference in the stability of the per-time-step alternating solution, as that is where $\epsilon \approx \epsilon_{\text{trunc}}$. The $\hat{\mathbf{x}}$-component of the truncated MOT-JVIE solution, $J^x(\mathbf{r},t)$, at $\mathbf{r} = (0.025,0.075,0.025)$ as a function of time is shown in Figure~\ref{fig:TruncationTest}. In Figure~\ref{fig:TruncationTest}, the per-time-step alternating solutions remains almost unaltered for $\epsilon_{\text{trunc}} \lessapprox 10^{-8}$. The first observed change is for $\epsilon_{\text{trunc}} \approx 10^{-7}$, which suggests that $\epsilon \approx 10^{-7}$. Then for $\epsilon_{\text{trunc}} \gtrapprox 10^{-6}$ the solution magnitude increases exponentially where the pertaining companion matrix eigenvalue increases with $\epsilon_{\text{trunc}}$.
\begin{Figure}
    \centering
    \includegraphics[width=\columnwidth]{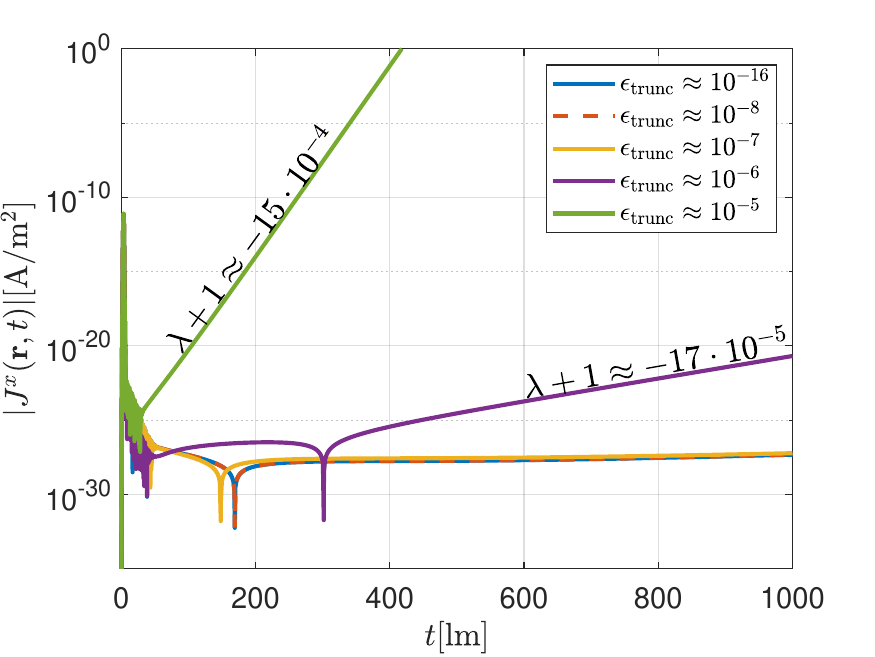}
    \captionof{figure}{
    The $\hat{\mathbf{x}}$-component of the truncated MOT-JVIE solution, $J^x(\mathbf{r},t)$, as a function of time sampled at $\mathbf{r} = (0.025,0.075,0.025)$ in the $0.2^3~\mathrm{m}^3$ dielectric cube with $\varepsilon_r = 2$. The discretization settings are $\Delta x = \Delta y = \Delta z = c_0 \Delta t = 0.01~\mathrm{m}$. The truncated MOT-JVIE has truncated interaction matrix values to deliberately introduce the truncation error $\epsilon_{\text{trunc}}$. The discretization settings are $\Delta x = \Delta y = \Delta z = c_0 \Delta t = 0.01~\mathrm{m}$. The companion matrix eigenvalues $\lambda$ associated to the nonphysical unstable solution are added to the respective lines.}
    \label{fig:TruncationTest}
\end{Figure}

Further minimization of the error $\epsilon$ to allow for a larger number of voxels is theoretically possible by improving the accuracy in e.g. the numerical integration over the test function in Equation~\eqref{eq:Z}, but that is difficult to obtain, because the produced magnetic fields are limited in smoothness~\cite{VanDiepen2024}. Therefore, we apply a simpler strategy, in which we enforce the positive-definiteness of $\hat{\mathbf{D}}_\ell$ via regularization and thus enforce stability of the solution pertaining the eigenvalues on the negative-real axis. However, regularization tends to result in a loss of accuracy in the solution. Improving the stability of MOT-schemes at the expense of accuracy has been proposed before, e.g. to the MOT-EFIE~\cite{Vechinski1992}. Here, we apply similar techniques based on the PDSA to improve stability while minimizing the loss in accuracy.

\subsection{Regularization}
Numerical experiments have shown that only $\mathbf{D}_\ell$ of the PDSA matrices $\mathbf{D}_n$ in Eq.~\eqref{eq:Dn} becomes indefinite when increasing the number of voxels in the MOT-JVIE and we conjecture this is due to finite precision errors in the elements of the interaction matrices. The matrix $\hat{\mathbf{D}}_\ell$~\eqref{eq:DellFinitePrecision} includes the finite precision errors in its formulation and, as discussed at the start of Section~\ref{sc:Stabilization}, its smallest eigenvalue is bounded by $\lambda^- - 3 M \epsilon$. To restore positive-definiteness of $\hat{\mathbf{D}}_\ell$, we apply regularization by adding a scaled identity matrix, i.e.
\begin{equation} \label{eq:Regularization}
    \hat{\mathbf{D}}^\delta_\ell = \hat{\mathbf{D}}_\ell + \delta \mathbf{I},
\end{equation}
where the lower bound of the smallest eigenvalue of $\hat{\mathbf{D}}^\delta_\ell$ is $\lambda^- - 3M\epsilon + \delta$. As discussed at the start of Section~\ref{sc:Stabilization}, the absolute error in the elements of the interaction matrix is around $\epsilon \approx 10^{-7}$. To compensate the finite precision error in the interaction matrix elements, we require the regularization parameter $\delta \gtrapprox 10^{-7} M$. Numerical experiments have confirmed that $\delta > 10^{-7} M$ is indeed sufficient to stabilize the MOT-JVIE.

We choose to implement the regularization in the MOT-scheme by replacing some of the interaction matrices $\mathbf{Z}_{n}$ by $\mathbf{Z}^\delta_{n} = \mathbf{Z}_{n} + \delta_{n} \mathbf{I}$. We choose $\delta_n$ to be only real numbers so the MOT-scheme remains real-valued. This regularization plays a roll in the sum of each row in Equation~\eqref{eq:MatrixEquation}, i.e.
\begin{equation}
    \sum_{n'= n-\ell} ^n \mathbf{Z}^\delta_{n-n'} \mathbf{J}_{n'} = \! \! \sum_{n'= n-\ell}^n \mathbf{Z}_{n-n'} \mathbf{J}_{n'} + \! \! \sum_{n'= n-\ell}^n \delta_{n-n'} \mathbf{J}_{n'},
\end{equation}
where the regularization with $\delta_n$ acts as a filter on the individual elements of $\mathbf{J}_n$ to increase the magnitude of some solution frequencies, $f$, while reducing that of others. The amount of regularization can be quantified as the normalized frequency magnitude response of a FIR-filter with real-valued coefficients $\delta_n$~\cite{Parks1987}, i.e.
\begin{equation} \label{eq:FIRmagnitude}
    |\sum_{n=0}^\ell e^{-\mathrm{j n \theta}} \delta_n|,
\end{equation}
where $\theta \in [0,\pi]$ is the normalized frequency, i.e. the solution frequency $f$ normalized to the time step size $\Delta t$ resulting in $\theta = 2 \pi \Delta t f$. The definition of the normalized frequency $\theta$ happens to coincide with the definition of the complex argument of the companion matrix eigenvalues~\cite{VanDiepen2024}. We refer to Equation~\eqref{eq:FIRmagnitude}
as the FIR-regularization magnitude.
To achieve the regularization of $\hat{\mathbf{D}}_\ell$ as in Equation~\eqref{eq:Regularization}, which concerns the normalized frequencies $\theta = \pi$, the values of $\delta_n$ should meet the criteria
\begin{equation} \label{eq:RegReq1}
    \sum_{n=0}^\ell (-1)^n \delta_n = \delta.
\end{equation}
This requirement only focuses on choosing $\delta_n$ such that it moves the companion matrix eigenvalues on the negative real axis from outside to inside the unit circle, i.e. the companion matrix eigenvalues with complex argument $\theta = \pi$. However, the other companion matrix eigenvalues, which represent the lower-frequency part of the solution, should preferably remain unaltered. Therefore, a second requirement on the regularization is the minimization of the FIR-regularization magnitude~\eqref{eq:FIRmagnitude}, i.e.
\begin{equation} \label{eq:RegReq2}
    \min_{\delta_n} |\sum_{n=0}^\ell e^{-\mathrm{j n \theta}} \delta_n| \text{ for } \theta \in [0,\pi) \text{ and } \delta_n \in \mathbb{R}.
\end{equation}
Although we try to minimize the complex-weighted sum over $\delta_n$ for all $\theta \neq \pi$ in the second requirement~\eqref{eq:RegReq2}, the FIR-regularization magnitude~\eqref{eq:FIRmagnitude} can still be close to $\delta$. In that case
there is a third requirement that the regularization should not move the companion matrix eigenvalues with complex argument $\theta \in [0,\pi)$ to the exterior of the unit circle. To meet these three requirements on the regularization, we base the values $\delta_n$ on real-valued low group-delay linear-phase high-pass FIR-filters~\cite{Parks1987}, because requirement~\eqref{eq:RegReq1} and \eqref{eq:RegReq2} are similar to design requirements in high-pass FIR filters. We empirically determined low group-delay linear phase is important to meet the third requirement. We refer to this type of regularization as FIR-regularization.

\subsection{Numerical results} \label{sc:StabilizationNumericalResults}
We consider three FIR-regularizations with increasing length in time, starting from a length of two time steps, up to four time steps. Consequently, the group phase delay increases as the filter length increases~\cite{Parks1987}. The two-step FIR-regularization (FIR2) is defined as,
\begin{equation} \label{eq:FIR2}
    \mathbf{Z}_0^\delta = \mathbf{Z}_0 + \frac{\delta}{2} \mathbf{I}, \ \mathbf{Z}_1^\delta = \mathbf{Z}_1 - \frac{\delta}{2} \mathbf{I},
\end{equation}
of which the regularization magnitude~\eqref{eq:FIRmagnitude} decreases linearly for $\theta \rightarrow 0$. The three-step FIR regularization (FIR3) is defined as,
\begin{equation} \label{eq:FIR3}
    \mathbf{Z}_0^\delta = \mathbf{Z}_0 + \frac{\delta}{4} \mathbf{I}, \ \mathbf{Z}_1^\delta = \mathbf{Z}_1 - \frac{\delta}{2} \mathbf{I}, \ \mathbf{Z}_2^\delta = \mathbf{Z}_2 + \frac{\delta}{4} \mathbf{I},
\end{equation}
of which the regularization magnitude~\eqref{eq:FIRmagnitude} decreases quadratically for $\theta \rightarrow 0$. The four-step FIR-regularization (FIR4) is defined as,
\begin{equation} \label{eq:FIR4}
\begin{split}
    \mathbf{Z}_0^\delta &= \mathbf{Z}_0 + \frac{\delta}{8} \mathbf{I}, \ \mathbf{Z}_1^\delta = \mathbf{Z}_1 - \frac{3\delta}{8} \mathbf{I}, \\ 
    \mathbf{Z}_2^\delta &= \mathbf{Z}_2 + \frac{3\delta}{8} \mathbf{I}, \ \mathbf{Z}_3^\delta = \mathbf{Z}_3 + \frac{\delta}{8} \mathbf{I}
\end{split}
\end{equation}
of which the regularization magnitude~\eqref{eq:FIRmagnitude} decreases cubically for $\theta \rightarrow 0$. 

We first apply these four FIR-regularzations to the experiment in Section~\ref{sc:AccelerationNumericalResults} for $\sqrt[3]{M} = 20$. We fix the regularization at $\delta = 0.1$ to accommodate for the largest discretization with $M = 10^6$. The $\hat{\mathbf{x}}$-component of the MOT-JVIE solution, $J^x(\mathbf{r},t)$, at $\mathbf{r} = (0.025,0.075,0.025)$ for the four different types of regularization are shown in Figure~\ref{fig:K20FIRComparision}.
\begin{Figure}
    \centering
    \includegraphics[width=\columnwidth]{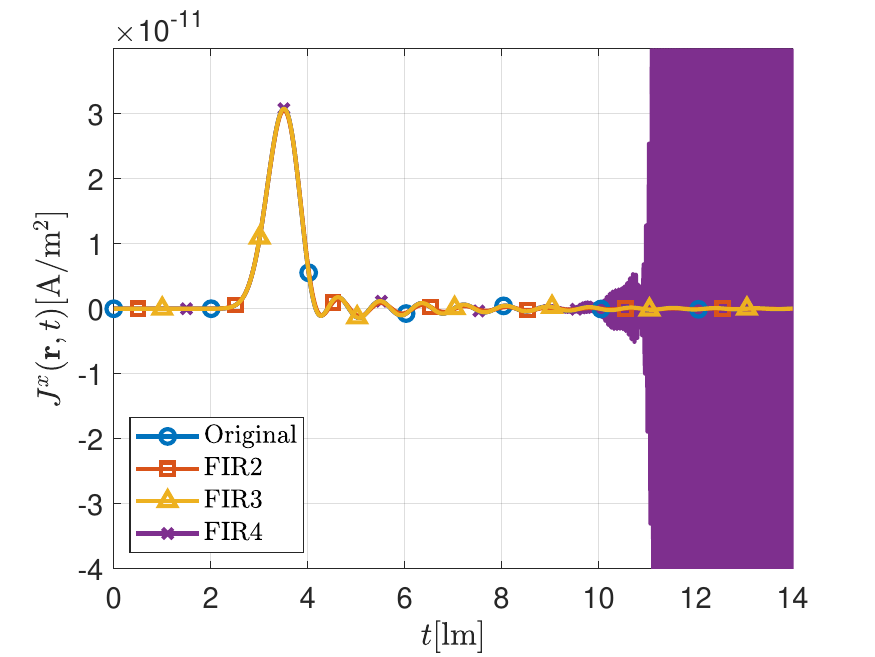}
    \captionof{figure}{The $\hat{\mathbf{x}}$-component of the $\delta = 0.1$ FIR2~\eqref{eq:FIR2}, FIR3~\eqref{eq:FIR3} and FIR4~\eqref{eq:FIR4}-regularized MOT-JVIE solution, $J^x(\mathbf{r},t)$, as a function of time sampled at $\mathbf{r} = (0.025,0.075,0.025)$ in the $0.2^3~\mathrm{m}^3$ dielectric cube with $\varepsilon_r = 12$. The discretization settings are $\Delta x = \Delta y = \Delta z = c_0 \Delta t = 0.2/\sqrt[3]{M}~\mathrm{m}$.}
    \label{fig:K20FIRComparision}
\end{Figure}
The result for FIR4-regularization illustrates what happens if the group delay of the regularization is too high, i.e. the regularization moves an companion matrix eigenvalue with complex argument $\theta \neq \pi$ to the exterior of the unit circle. The absolute difference between the solutions of the FIR-regularized MOT-JVIE and the original MOT-JVIE normalized to the value $J^x = 3\cdot 10^{-11} [\mathrm{A/m^2}]$, i.e. approximately the peak solution, is shown in Figure~\ref{fig:K20FIRDifference}.
\begin{Figure}
    \centering
    \includegraphics[width=\columnwidth]{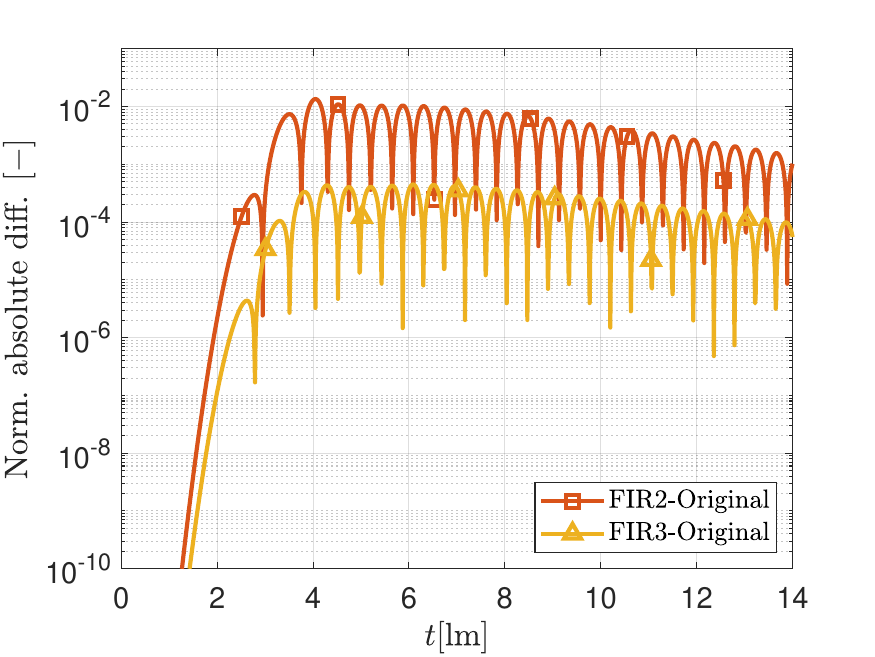}
    \captionof{figure}{The absolute difference between the FIR-regularized and original MOT-JVIE solution shown in Figure~\ref{fig:K20FIRComparision} normalized to the value $J^x = 3\cdot 10^{-11} [\mathrm{A}/\mathrm{m}^2]$.}
    \label{fig:K20FIRDifference}
\end{Figure}
There we observe that the faster decay of the regularization with respect to $\theta \rightarrow 0$, i.e. FIR2~\eqref{eq:FIR2} decays quadratically and FIR3~\eqref{eq:FIR3} cubically, diminishes the error introduced by the regularization.

To further test the $\delta = 0.1$ FIR3-regularization~\eqref{eq:FIR3}, we also apply it to the other discretizations in the experiment in Section~\ref{sc:AccelerationNumericalResults}, i.e. $\sqrt[3]{M} \in [36,52,84,100]$. The $\hat{\mathbf{x}}$-component of the contrast current density at $\mathbf{r} = (0.025,0.075,0.025)$ is shown in Figure~\ref{fig:Cube_epsr12_StabilityComparison}.
\begin{Figure}
    \centering
    \includegraphics[width=\columnwidth]{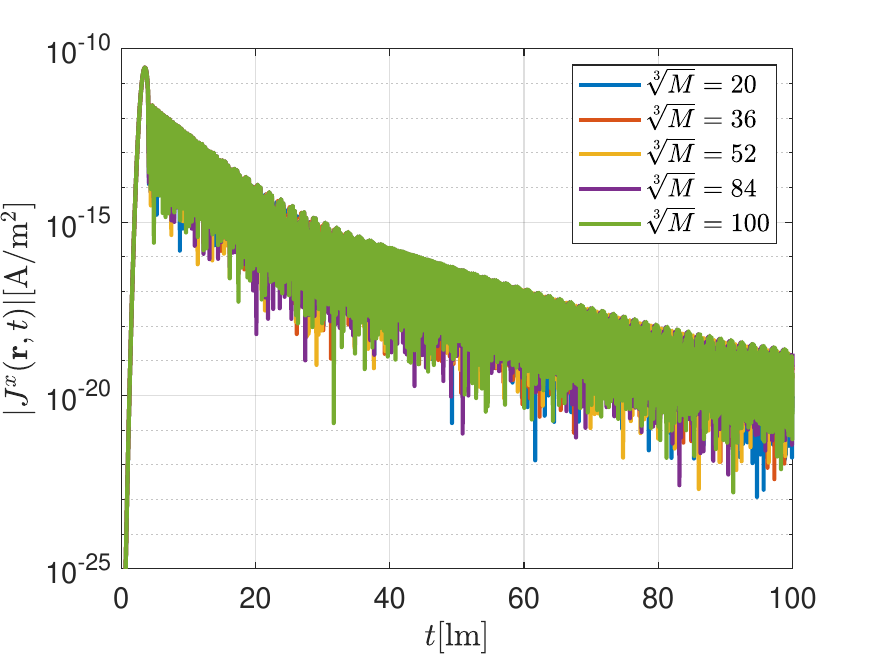}
    \captionof{figure}{The $\hat{\mathbf{x}}$-component of the $\delta = 0.1$ FIR3-regularized MOT-JVIE solution, $J^x(\mathbf{r},t)$, as a function of time sampled at $\mathbf{r} = (0.025,0.075,0.025)$ in the $0.2^3~\mathrm{m}^3$ dielectric cube with $\varepsilon_r = 12$. The discretization settings are $\Delta x = \Delta y = \Delta z = c_0 \Delta t = 0.2/\sqrt[3]{M}~\mathrm{m}$.}
    \label{fig:Cube_epsr12_StabilityComparison}
\end{Figure}
The FIR3-regularized MOT-JVIE remains stable in the simulation time span, unlike the MOT-JVIE solution without regularization in Figure~\ref{fig:Cube_epsr12_InstabilityComparison}. The average absolute difference between the $\delta = 0.1$ FIR3-regularized MOT-JVIE and the original MOT-JVIE normalized to the value $J^x = 3 \cdot 10^{-11}$ in the $64$-sample locations is shown in Figure~\ref{fig:FIR3AvgL2Difference} for $\sqrt[3]{M} \in [20,52,100]$. This value increases for $\sqrt[3]{M} = 100$ at $t \geq 12~\mathrm{lm}$, because the nonphysical unstable original MOT-JVIE solution becomes dominant, see Figure~\ref{fig:Cube_epsr12_InstabilityComparison}. Still, the error introduced by the regularization decreases approximately as $1/8$ when approximately doubling $\sqrt[3]{M}$. The decrease of the error is due the halving of the time step between simulations as $\Delta t = 0.2/\sqrt[3]{M}$. Subsequently, the complex argument $\theta$ of all companion matrix eigenvalues halves. The FIR3-regularization magnitude~\eqref{eq:FIR3} decreases cubically with $\theta$, hence the factor $1/8$.

\begin{Figure}
    \centering
    \includegraphics[width=\columnwidth]{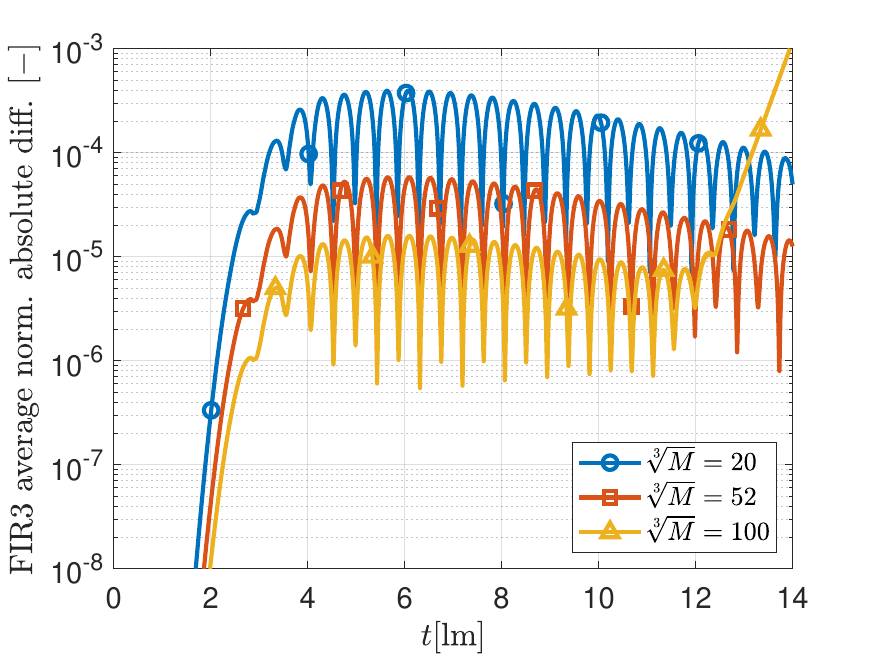}
    \captionof{figure}{The absolute difference between the $\delta = 0.1$ FIR3-regularized and original MOT-JVIE solution normalized to the value $J^x = 3\cdot 10^{-11} [\mathrm{A}/\mathrm{m}^2]$ averaged over the $64$ sample points.}
    \label{fig:FIR3AvgL2Difference}
\end{Figure}

\section{FFT-Accelerated stabilized MOT-JVIE} \label{sc:AccRegMOTJVIE}
We have FFT-accelerated and stabilized the MOT-JVIE presented in~\cite{VanDiepen2024} as explained in Section~\ref{sc:SpatialTemporalAcceleration} and Section~\ref{sc:Stabilization}, respectively. We will now test the capibilities of the FIR3-regularized MOT-JVIE.

\subsection{Cube} \label{sc:Cube}
We start by evaluating the accuracy of the numerical experiment with $\sqrt[3]{M} = 100$ in Section~\ref{sc:StabilizationNumericalResults}, i.e. the contrast current density in a $0.2^3 ~\mathrm{m}^3$ cubic domain with $\varepsilon_r = 12$ induced by a Gaussian $\hat{\mathbf{x}}$-polarized plane wave travelling in the negative $\hat{\mathbf{z}}$-direction measured at $64$ sample locations. We compare the $\hat{\boldsymbol{\alpha}}$-component of the frequency magnitude response, $|H^\alpha|$, based on the FIR3-regularized MOT-JVIE solution to the frequency magnitude response computed with CST Studio Suite 2023~\cite{CST2023}, $|H_{\mathrm{CST}}^\alpha|$, using its combined-field integral equation frequency-domain solver. The $\hat{\boldsymbol{\alpha}}$-component of the frequency magnitude response is defined as 
\begin{equation} \label{eq:Halpha}
    |H^\alpha|(\mathbf{r},f) = \frac{|\hat{\boldsymbol{\alpha}} \cdot \mathbf{j}_\varepsilon(\mathbf{r},f)|}{\varepsilon_0(\varepsilon_r-1)|e^i(f)|},
\end{equation}
where $\mathbf{j}_\varepsilon$ is the frequency-domain counterpart of the $\mathbf{J}_\varepsilon$~\eqref{eq:ContrastCurrentDensity} and the $|e^i|$ is the magnitude of the Gaussian plane wave in the frequency domain. The computation of both $\mathbf{j}_\varepsilon$ and $|e^i|$ are explained in~\cite{VanDiepen2024}, however, we analyze this problem at higher frequencies with corresponding stronger resonances compared to~\cite{VanDiepen2024}, i.e. the resulting frequency magnitude response has much narrower and has higher peaks, which correspond to time-domain solutions that decay slowly over time as observed in Figure~\ref{fig:Cube_epsr12_StabilityComparison}. The limited simulation time, in Figure~\ref{fig:Cube_epsr12_StabilityComparison} up to $T = 100~\mathrm{lm}$, is therefore an abrupt truncation of the time-domain solution that reduces the accuracy of the computed magnitude response~\cite{Bloomfield2000}. Instead of running the simulation for a longer time, we taper the time-domain solution. The tapering of time domain signals has been analyzed in~\cite{Bloomfield2000} and we adopt one of the recommendations: a tapered cosine window on the last $20~\%$ of the time-domain samples.

The three Cartesian components of the frequency magnitude response, $|H^\alpha|(\mathbf{r},f)$ for $\alpha \in [x,y,z]$, at $\mathbf{r} = (0.025,0.075,0.025)$ as a function of frequency are shown in Figure~\ref{fig:K100Halpha}. The CST reference based on the combined field integral equation~\cite{CST2023}, was unable to produce results for frequencies lower than $100~\mathrm{MHz}$. The bandwidth of the Gaussian plane wave is limited by the double-precision-arithmetic noise floor, which we approach here at $f \approx 900~\mathrm{MHz}$. Therefore, we truncate the result at $f = 900~\mathrm{MHz}$ in Figure~\ref{fig:K100Halpha}.
\begin{Figure}
    \centering
    \includegraphics[width=\columnwidth]{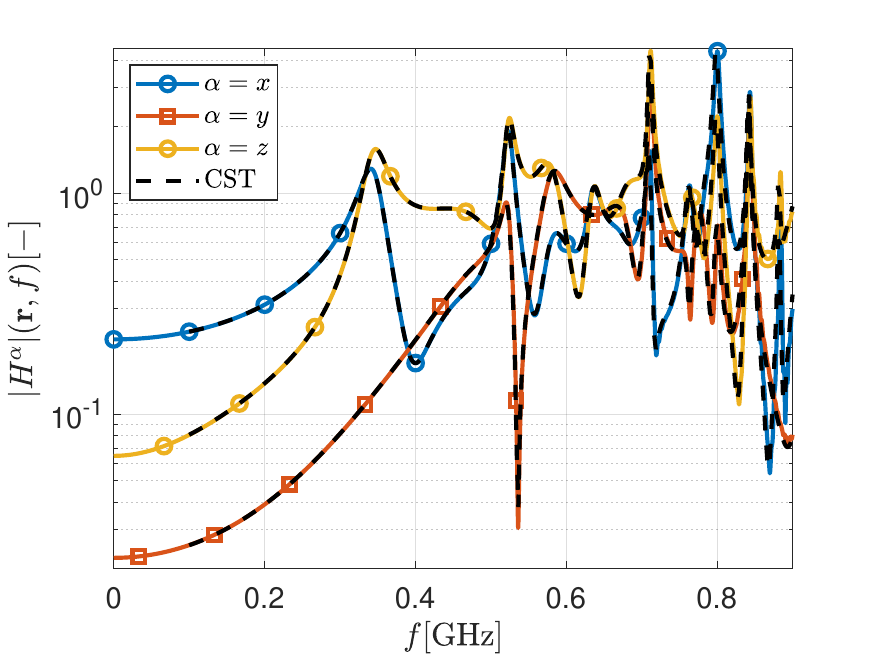}
    \captionof{figure}{The frequency magnitude response $\hat{\boldsymbol{\alpha}}$-component of the $\delta = 0.1$ FIR3-regularized MOT-JVIE solution, $|H^\alpha(\mathbf{r}, f)|$, as a function of frequency sampled at $\mathbf{r} = (0.025,0.075,0.025)$ in the $0.2^3~\mathrm{m}^3$ dielectric cube with $\varepsilon_r = 12$. The discretization settings are $\Delta x = \Delta y = \Delta z = c_0 \Delta t = 0.002~\mathrm{m}$.}
    \label{fig:K100Halpha}
\end{Figure}
In Figure~\ref{fig:K100Halpha}, we observe the overlap between the MOT-JVIE result and that of CST. To quantify accuracy, we compute the $L^2$-relative error defined as
\begin{equation} \label{eq:L2error}
    L^2(f) = \sqrt{\frac{\sum (|H^\alpha(\mathbf{r},f)|-|H_{\mathrm{CST}}^\alpha(\mathbf{r},f)|)^2}{\sum (|H_{\mathrm{CST}}^\alpha(\mathbf{r},f)|)^2}},
\end{equation}
where the summation is over the $x$, $y$ and $z$ components and the aforementioned $64$ sample points. The $L^2$-relative error as a function of frequency is shown in Figure~\ref{fig:L2RelativeNorm}.
\begin{Figure}
    \centering
    \includegraphics[width=\columnwidth]{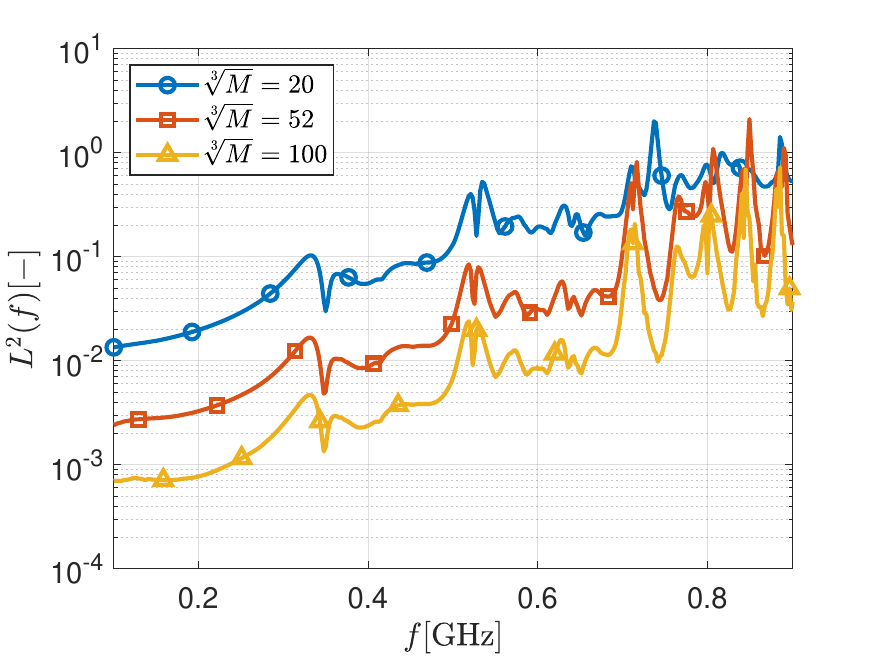}
    \captionof{figure}{The $L^2$-relative error as defined in~\eqref{eq:L2error} based on the $\delta = 0.1$ FIR3-regularized MOT-JVIE solution.}
    \label{fig:L2RelativeNorm}
\end{Figure}
In Figure~\ref{fig:L2RelativeNorm} we observe local peaks in the accuracy and an overall increasing trend. These local peaks occur near the resonances in the solution, where a small mismatch in the peak frequency results in large relative errors. The overall trend is that solution accuracy improves when the number of voxels is increased, which corresponds to an increase in the number of voxels per wavelength. This is in line with what we expect, i.e. the limiting factor in accuracy is still the number of voxels per wavelength~\cite{VanDiepen2024}. 

\subsection{Sphere} \label{sc:Sphere}
The second test consists of the FIR3-regularized MOT-JVIE applied to a sphere with high permittivity $\varepsilon_r = 100$. This type of scatterer is difficult for the MOT-JVIE for two reasons. The first reason is that the curvature of the outer surface results in a stair-casing error in the discretization, as illustrated in Figure~\ref{fig:Voxelizev3}. The second reason is that the sphere has even stronger resonances, i.e. higher and  narrower peaks in the frequency magnitude response, than those observed for the $\varepsilon_r = 12$ cube in Section~\ref{sc:Cube}. We observed that these strong resonances are the most difficult to match in frequency with the FIR3-regularized MOT-JVIE. Therefore we want to further evaluate the performance of the MOT-JVIE.

The sphere has a diameter of $0.2~\mathrm{m}$ and is centered at $\mathbf{r} = (0.1,0.1,0.1)$. The contrast current density inside the sphere is induced by the Gaussian plane wave in~\eqref{eq:GaussPlaneWave} with $\sigma = 2~\mathrm{lm}$, $t_0 = 3.42~\mathrm{lm}$, $E_0 = 1~\mathrm{V/m}$, $\hat{\mathbf{k}} = (-1/2,-1/2, -\sqrt{2}/2)$ and $\hat{\mathbf{p}} = (1/2,1/2, -\sqrt{2}/2)$. We voxelize the sphere as explained in Section~\ref{sc:Voxelization}, i.e. the sphere is enclosed in a $0.2^3~\mathrm{m}^3$ box and this box is evenly divided in $105^3$ voxels with edge lengths $\Delta x = \Delta y = \Delta z = 0.2/105~\mathrm{m}$. The discrete time step size is set to $ \Delta t = 0.2/105~\mathrm{lm}$. We choose $\delta = 105^3 \cdot 10^{-7}$ for the FIR3-regularization~\eqref{eq:FIR3} to maintain stability, as explained in Section~\ref{sc:Stabilization}.  
\begin{Figure}
    \centering
    \includegraphics[width=\columnwidth]{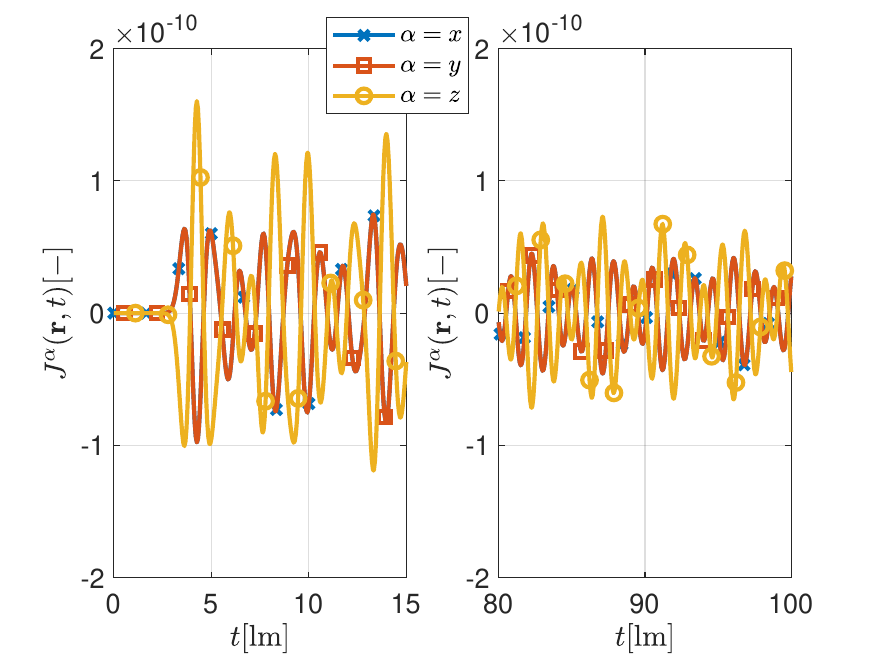}
    \captionof{figure}{The $\hat{\boldsymbol{\alpha}}$-component of the $\delta = 105^3 \cdot 10^{-7}$ FIR3-regularized MOT-JVIE solution, $J^\alpha(\mathbf{r},t)$, as a function of time sampled at $\mathbf{r} = (0.1 + 0.1/105,0.1 + 0.1/105,0.1 + 0.1/105)$ in the $0.2~\mathrm{m}$ diameter sphere with $\varepsilon_r = 100$. The discretization settings are $\Delta x = \Delta y = \Delta z = c_0 \Delta t = 0.2/105~\mathrm{m}$.}
    \label{fig:SphereK105J}
\end{Figure}

From the contrast current density solution in Figure~\ref{fig:SphereK105J}, we compute the frequency magnitude response defined as
\begin{equation} \label{eq:H}
    |H|(\mathbf{r},f) = \! \sqrt{|H^x|(\mathbf{r},f)^2+|H^y|(\mathbf{r},f)^2+|H^z|(\mathbf{r},f)^2},
\end{equation}
where $|H^x|$, $|H^y|$ and $|H^z|$ are the three Cartesian components of the frequency magnitude response in~\eqref{eq:Halpha}. As the resonances in the $\varepsilon_r = 100$ sphere are even stronger than those of the $\varepsilon_r = 12$ cube in Section~\ref{sc:Cube}, we need to extend the tapering of the solution from the last $20~\%$ of the solution in Section~\ref{sc:Cube} to almost the entire length solution. This significantly reduces the accuracy of the computation as it flattens the resonant peaks~\cite{Bloomfield2000}. Therefore, we use an alternative technique, i.e. vector-fitting of a transfer function on the time-domain sequence, from which we subsequently determine the frequency magnitude response. The {\tt tfest}-function in MATLAB~\cite{MATLAB2023} performs this vector-fitting and yields a transfer function with $40$ poles and $39$ zeros to the time-domain data in Figure~\ref{fig:SphereK105J}. The pertaining frequency magnitude response together with the one obtained from the Mie series is shown in Figure~\ref{fig:SphereK105H}.
\begin{Figure}
    \centering
    \includegraphics[width = \columnwidth]{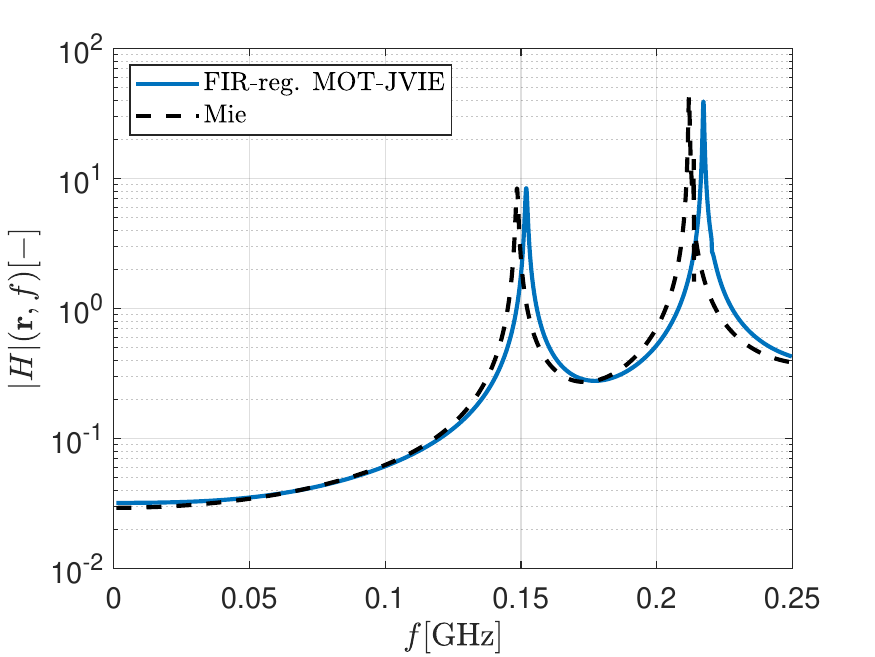}
    \captionof{figure}{The frequency magnitude response, $|H|(\mathbf{r}, f)$~\eqref{eq:H}, as a function of frequency based on the Cartesian components of the FIR3-regularized MOT-JVIE contrast current density solution in Figure~\ref{fig:SphereK105J} obtained via vector fitting, and the one obtained via the Mie-series.}
    \label{fig:SphereK105H}
\end{Figure}
The first resonance peak of the FIR3-regularized solution in Figure~\ref{fig:SphereK105H} has a $2.3\%$ shift in frequency and a $2.6\%$ relative error in magnitude, compared to the Mie-series solution. The second resonance peak of the FIR3-regularized solution in Figure~\ref{fig:SphereK105H} has a $2.6\%$ frequency shift and a $8.6\%$ relative error in magnitude, compared to the Mie-series solution.

\section{Conclusion}~\label{sc:Conclusion}
We focused on two parts in this work concerning the marching-on-in-time contrast current density volume integral equation (MOT-JVIE), i.e. the fast-Fourier-transform (FFT) acceleration and the stabilization through regularization of the MOT-JVIE. We reviewed the existing FFT-acceleration of time-domain surface integral equations, where we observed that both a hierarchical division in space and hierarchical division in time result in the same $\mathcal{O}(M \log^2 M)$-scaling in computation time, where $M$ is the number of spatial unknowns. We extended the hierarchical division in space that was already available in two dimensions in literature to the three-dimensional MOT-JVIE and demonstrated the $\mathcal{O}(M \log^2 M)$-scaling by numerical experiments. We employed the positive definite stability analysis (PDSA) to analyze an observed instability in the MOT-JVIE related to an increase in the number of voxels. A link between stability, finite precision in the matrix elements, and the number of voxels was observed and we concluded that the number of voxels for a stable MOT-JVIE discretization is restricted by the finite precision of the matrix elements. The analysis with the PDSA showed that stability can be enforced through regularization of the MOT-JVIE, at the cost of accuracy. We introduced FIR-regularization based on low group-delay linear-phase high-pass FIR-filters, to minimize the impact on the accuracy, and we illustrated the impact with numerical experiments. The capabilities of the FFT-accelerated FIR-regularized MOT-JVIE were illustrated by computing the time-domain results for a high-permittivity cube and sphere discretized with approximately a million voxels. These results were converted from the time domain to the frequency domain and compared to those obtained via a commercial combined-field integral equation solver and the Mie series, respectively.

\appendix
\section{Toeplitz FFT-acceleration} \label{ap:ToeplitzFFTAcceleration}

Consider the $n_{\text{row}} \times n_{\text{col}}$ Toeplitz matrix $\mathbf{A}$ where $n_{\text{row}} \geq n_{\text{col}}$
\begin{equation}
	\mathbf{A} = \begin{bmatrix}
        \mathbf{A}_0 & \mathbf{A}_{-1} & \cdots & \mathbf{A}_{-n_{\text{col}}+1} \\
        \mathbf{A}_{1} & \ddots & \ddots & \vdots \\
        \vdots & \ddots & \ddots & \mathbf{A}_{-1} \\
        \mathbf{A}_{n_{\text{row}}-n_{\text{col}}} &  & \ddots & \mathbf{A}_0 \\
        \vdots & \ddots& &  \mathbf{A}_{1} \\
        \vdots & & \ddots & \vdots \\
        \mathbf{A}_{n_{\text{row}}-1}& \cdots & \cdots & \mathbf{A}_{n_{\text{row}}-n_{\text{col}}} \\
        
	\end{bmatrix}.
\end{equation}
The following analysis also holds if the elements of $\mathbf{A}$ are themselves square matrices, i.e. when $\mathbf{A}$ is a block-Toeplitz matrix. To keep the explanation concise, we only perform the analysis for $\mathbf{A}$ with scalar matrix elements. The case where $n_{\text{row}} = n_{\text{col}}$, i.e. $\mathbf{A}$ is a square Toeplitz matrix, is well documented~\cite{Golub1996}, but we extend this analysis to non-square Toeplitz matrices.

The Toeplitz matrix $\mathbf{A}$ is uniquely defined by the elements in the first row and column. A circulant matrix $\mathbf{C}$ is uniquely defined by the elements in the first column. If we define the first column of $\mathbf{C}$ as
\begin{equation}
    \mathbf{C}_{\text{col}} = [\mathbf{A}_0;\ldots;\mathbf{A}_{n_{\text{row}}-1};\mathbf{A}_{-n_{\text{col}}+1};\ldots;\mathbf{A}_{-1}],
\end{equation}
where $;$ is the separation between column elements, then the top-left block of $\mathbf{C}$ is equivalent to $\mathbf{A}$. The resulting circulant matrix $\mathbf{C}$ is square with a dimension $n_{\text{col}}+n_{\text{row}}-1$. To compute the matrix vector product of  for a known vector $\mathbf{x}$, i.e. $\mathbf{b}= \mathbf{A}\mathbf{x}$, one can replace this by
\begin{equation} \label{eq:ToeplitzToCirculant}
    \mathbf{b} = \mathbf{A}\mathbf{x} \rightarrow
    \mathbf{C}
    \begin{bmatrix}
        \mathbf{x}\\
        \mathbf{0}
    \end{bmatrix}
    =
    \begin{bmatrix}
        \mathbf{A} & \cdot \\
        \cdot & \cdot
    \end{bmatrix}
    \begin{bmatrix}
        \mathbf{x}\\
        \mathbf{0}
    \end{bmatrix} = \begin{bmatrix}
        \mathbf{b}\\
        \mathbf{c}
    \end{bmatrix},
\end{equation}
where $[\mathbf{x};\mathbf{0}]$ is the zero-padding of $\mathbf{x}$ from a length $n_{\text{col}}$ to $N$ and $\mathbf{c}$ is the side effect of replacing $\mathbf{A}$ by $\mathbf{C}$. So, we have shown that an $n_{\text{row}} \times n_{\text{col}}$ Toeplitz matrix with $n_{\text{row}} \geq n_{\text{col}}$ can be viewed as the top-left truncation of a circulant matrix of dimension $N\times N$ with $N=n_{\text{row}} + n_{\text{col}}-1$. These steps can be repeated for a $n_{\text{row}} \times n_{\text{col}}$ Toeplitz matrix with $n_{\text{row}} \leq n_{\text{col}}$, which also results in a circulant matrix of the same dimension. 

The number of operations required to compute $\mathbf{A} \mathbf{x}$ for a known vector $\mathbf{x}$ scales with the number of elements in $\mathbf{A}$, i.e. $\mathcal{O}(n_{\text{row}} n_{\text{col}})$. Even though the circulant matrix $\mathbf{C}$ is larger than $\mathbf{A}$, the computation of $\mathbf{C} [\mathbf{x};\mathbf{0}]$ can be faster because the matrix is circular. By applying a fast Fourier transform (FFT) and its inverse (IFFT)~\cite{Golub1996}, we can rewrite the matrix vector product involving the circulant matrix $\mathbf{C}$ in Equation~\eqref{eq:ToeplitzToCirculant} as
\begin{equation}
    \begin{bmatrix}
        \mathbf{b}\\
        \mathbf{c}
    \end{bmatrix} = \mathrm{IFFT}\left(\mathrm{FFT}(\mathbf{C}_{\text{col}}) \odot \mathrm{FFT}\left(\begin{bmatrix}
        \mathbf{x}\\
        \mathbf{0}
    \end{bmatrix}\right) \right),
\end{equation}
where $\odot$ is the point-wise multiplication of two column vectors and $\mathrm{FFT}(\cdot)$ and $\mathrm{IFFT}(\cdot)$ are the FFT and its inverse operators, respectively, which definitions can be found in~\cite{Golub1996} and implementations in~\cite{Frigo2005}. The FFT and IFFT operations scale as $\mathcal{O}((n_{\text{col}}+n_{\text{row}})\log(n_{\text{col}}+n_{\text{row}}))$ and the point-wise multiplication as $\mathcal{O}(n_{\text{col}}+n_{\text{row}})$~\cite{Frigo2005}, thus the FFT dominates this operation. The FFT-acceleration, i.e. replacing the Toeplitz matrix vector product by a circulant matrix vector product, scales as 
\begin{equation} \label{eq:FFTAcceleration}
    \mathcal{O}\left(\frac{n_{\text{row}} n_{\text{col}}}{(n_{\text{row}}+n_{\text{col}})\log(n_{\text{row}}+n_{\text{col}})}\right).
\end{equation}
The above equation illustrates that the FFT-acceleration actually accelerates as long as $\log n_{\text{col}} < n_{\text{row}}$ or $\log n_{\text{row}} < n_{\text{col}}$, and reaches it optimum  roughly when $n_{\text{row}} \approx n _{\text{col}}$, i.e. when the Toeplitz matrix is approximately square. If the $N \approx n_{\text{row}} \approx n _{\text{col}}$, the FFT-acceleration scaling simplifies to $\mathcal{O}\left(\frac{N}{\log N}\right)$, which illustrates that FFT-acceleration scales proportional to $N$ for an (approximately) square $N\times N$ Toeplitz matrix.

In the case $\mathbf{A}$ is a block-Toeplitz matrix, the acceleration via FFTs is still possible in the direction of the Toeplitz structure. If the block matrices are also block-Toeplitz, the procedure repeats itself for each level of the block-Topelitz structure, where one applies multi-dimensional FFTs where the dimension corresponds to the number of levels~\cite{Catedra1989,Zwamborn1991}. The complexity of matrix vector product scales then with the total size of the multi-level block-Toeplitz matrix.

\nocite{*}
\bibliographystyle{ieeetr}
\bibliography{Bibliography}

\end{multicols}
\end{document}